\documentclass[a4paper]{article}

\usepackage{epsfig}
\usepackage{subfigure}

\usepackage[latin1]{inputenc}
\usepackage[english]{babel}

\usepackage{longtable}
\usepackage{booktabs}

\newif\ifdeuxcols
\makeatletter \if@twocolumn
  \deuxcolstrue
\else
  \deuxcolsfalse
\fi \makeatother

\makeatletter
\newcommand{\purpose}[1]{\def\@purpose{#1}}
\purpose{}
\def\@oddfootlefttext{
  \ifx\@purpose\@empty
    Preprint submitted to \ifx\@journal\@empty Elsevier\else\@journal\fi
  \else\@purpose\fi}

\def\ps@pprintTitle{%
     \let\@oddhead\@empty
     \let\@evenhead\@empty
     \def\@oddfoot{\footnotesize\itshape
       \@oddfootlefttext\hfill\today}%
     \let\@evenfoot\@oddfoot}
\makeatother

\usepackage{amsmath} 
\usepackage{amssymb}  
\usepackage{epsfig}
\usepackage{graphicx}
\usepackage[hang]{caption2}
\newtheorem{remark}{Remark}[section]
\newtheorem{example}{Example}[section]

\usepackage{natbib}

\begin{document}

\title{\Huge \bf Model-free control}

\author{Michel FLIESS$^{a,c}$, C\'{e}dric JOIN$^{b,c,d}$ \\
  \\
  \multicolumn{1}{p{.9\textwidth}}{
  \centering\emph{
  $^{a}${\em {LIX (CNRS, UMR 7161)\\ \'Ecole polytechnique, 91128 Palaiseau, France \\ {\tt Michel.Fliess@polytechnique.edu}}}\\  \vspace{0.3cm}
  $^{b}${\em CRAN (CNRS, UMR 7039)\\ Universit\'{e} de Lorraine, BP 239 \\ 54506 Vand{\oe}uvre-l\`{e}s-Nancy, France \\{\tt cedric.join@univ-lorraine.fr}}\\ \vspace{0.3cm}
  $^{c}${\em AL.I.E.N.\\ {\em (ALgèbre pour Identification et Estimation Numériques)} \\ 24-30 rue Lionnois, BP 60120, 54003 Nancy, France\\ {\tt \{michel.fliess, cedric.join\}@alien-sas.com}}\\  \vspace{0.3cm}
  $^{d}${\em Projet NON-A, INRIA Lille -- Nord-Europe, France}
  }}}

%
%

\date{}

\maketitle

\begin{abstract}
``Model-free control'' and the corresponding ``intelligent'' PID
controllers (iPIDs), which already had many successful concrete
applications, are presented here for the first time in an unified
manner, where the new advances are taken into account. The basics
of model-free control is now employing some old functional analysis
and some elementary differential algebra. The estimation techniques
become quite straightforward via a recent online parameter
identification approach. The importance of iPIs and especially of iPs is deduced from the
presence of friction.  The strange industrial ubiquity of classic
PID's and the great difficulty for tuning them in complex situations is deduced, via an elementary sampling, from their connections with iPIDs. Several
numerical simulations are presented which include some
infinite-dimensional systems. They demonstrate not only the power of our intelligent controllers but also the great simplicity for tuning them.\\
\\

\end{abstract}

\newpage
\section{Introduction}\label{introduction}
Although \emph{model-free control} was introduced only a few years ago (\cite{esta,malo,marseille}),
there is already a quite impressive list of successful concrete
applications in most diverse fields, ranging from intelligent
transportation systems to energy management
(\cite{sofia,acc,cifa10,buda,cifa12,it13,nolcos,brest,psa,edf0,edf,lafont,michel,mil,mex,
vil1,vil2,vil3,vil4,mounier}). Most of those references lead
to practical implementations. Some of them are related to patents.

\begin{remark}
The wording \textit{model-free control} is of course not new in the literature, where it
has already been employed by a number of authors. The corresponding literature is huge: see,
\textit{e.g.}, \cite{kadri,chang,hahn,hong,keel,kill,malis,santos,spall,belg,sya,xu}. The corresponding settings are quite varied. They range from ``classic'' PIDs to robust and adaptative control via techniques stemming 
from, \textit{e.g.}, neural nets, fuzzy systems, and soft computing.  To
the best of our understanding, those approaches are rather
far from what we are developing here. Let us emphasize however Remark \ref{remote} for a comment on some works which are perhaps closer. See also Remark \ref{global}.
\end{remark}

Let us now summarize some of the main theoretical ideas which are shaping our model-free control. We restrict ourselves for simplicity's sake 
to systems with a single control variable $u$ and a single output variable $y$.
The unknown ``complex'' mathematical model is replaced by an \emph{ultra-local
model}
\begin{equation}
\boxed{y^{(\nu)} = F + \alpha u} \label{ultralocal}
\end{equation}
\begin{enumerate}
\item $y^{(\nu)}$ is the derivative of order $\nu \geq 1$ of $y$. The integer $\nu$
is selected by the practitioner. The existing
examples show that $\nu$ may always be chosen quite low,
\textit{i.e.}, $1$, or, only seldom, $2$. See Section
\ref{ord1} for an explanation.
\item $\alpha \in \mathbb{R}$ is a non-physical constant parameter.
It is chosen by the practitioner such that $\alpha u$ and
$y^{(\nu)}$ are of the same magnitude. It should be therefore clear that its numerical value, which is obtained by trials and errors, is not \textit{a priori} precisely defined. Let us stress moreover 
that controlling industrial plants has always been achieved by collaborating with engineers who know the system behaviour well.
\item $F$, which is continuously updated,  subsumes the poorly known parts of the plant as well as of the
various possible disturbances, without the need to make any
distinction between them.
\item For its estimation, $F$ is approximated by a piecewise constant function. Then the algebraic identification techniques due to \cite{sira1,sira2} are applied to the equation
\begin{equation}\label{identif}
y^{(\nu)} = \phi + \alpha u
\end{equation}
where $\phi$ is an unknown constant parameter. The estimation
\begin{itemize}
\item necessitates only a quite short time lapse,
\item is expressed via algebraic formulae which contain low-pass filters like iterated time integrals,
\item is robust with respect to quite strong noise corruption, according to the new setting of noises via \emph{quick fluctuations} (\cite{bruit}). 
\end{itemize}
\end{enumerate}
\begin{remark}
The following comparison with computer graphics might be enlightening. Reproducing on a screen a complex plane curve is not achieved via the equations defining that curve but by approximating
it with short straight line segments. Equation \eqref{ultralocal} might be viewed as a kind of analogue of such a short segment.
\end{remark}
\begin{remark}\label{global}
Our terminology \emph{model-free control} is best explained by the ultra-local Equation \eqref{ultralocal} which implies that the need of any ``good'' and ``global'' modeling is abandoned.
\end{remark}
 Assume that $\nu = 2$ in Equation \eqref{ultralocal}:
\begin{equation}
\ddot{y} = F + \alpha u \label{2}
\end{equation}
Close the loop via the \emph{intelligent
proportional-integral-derivative controller}, or \emph{iPID},
\begin{equation}\label{ipid}
\boxed{u = - \frac{F - \ddot{y}^\ast + K_P e + K_I \int e + K_D
\dot{e}}{\alpha}}
\end{equation}
where 
\begin{itemize}
\item $y^\ast$ is the reference trajectory,
\item $e = y - y^\ast$ is the tracking error,
\item $K_P$, $K_I$, $K_D$ are the usual tuning gains.
\end{itemize}
Combining
Equations \eqref{2} and \eqref{ipid} yields
\begin{equation}\label{stabil}
\ddot{e} + K_D \dot{e} + K_P e + K_I \int e = 0
\end{equation}
Note that $F$ does not appear anymore in Equation \eqref{stabil}, {\it i.e.,} the unknown parts and disturbances
of the plant vanish. We are therefore left with a linear differential equation with constant coefficients of order $3$. The tuning of $K_P$, $K_I$, $K_D$ becomes
therefore straightforward for obtaining a ``good'' tracking of $y^\ast$. This is a major benefit when
compared to the tuning of ``classic'' PIDs.
\begin{remark}
\emph{Intelligent PID controllers} may already be found in the literature but with a different meaning (see, \textit{e.g.}, \cite{astrom2}).
\end{remark}
\begin{remark}\label{remote}
See, \textit{e.g.}, \cite{chang-jung,han,youcef,zheng} for some remote analogy with our calculations. 
Those references assume however that the system order is finite and moreover known.
\end{remark}

Our paper is organized as follows. The general principles of
model-free control and of the corresponding intelligent PIDs are
presented in Section \ref{mfc}. The online estimation of the crucial
term $F$ is discussed in Section \ref{2}. Section \ref{ord1}
explains why the existence of frictions permits to restrict our
intelligent PIDs to intelligent proportional or to intelligent
proportional-integral correctors. The numerical
simulations in Section \ref{ex} examine the
following case-studies:
\begin{itemize}
\item A part of the unknown system may be nevertheless known. If it
happens to be {\em flat} (\cite{flmr}, and \cite{levine,sira}), it will greatly
facilitate the choice of the reference trajectory and of the
corresponding nominal control variable.
\item Standard modifications including aging and an actuator fault keep the
performances, with no damaging, of our model-free control synthesis without the need
of any new calibration.
\item An academic nonlinear case-study demonstrates that a single model-free
control is sufficient whereas many classic PIDs may be necessary in
the usual PID setting.
\item Two examples of infinite-dimensional systems demonstrate that our model-free control
provides excellent results without any further ado:
\begin{itemize}
\item a system with varying delays,
\item a one-dimensional semi-linear heat equation, which is borrowed
from \cite{Coron}.
\end{itemize}
\item A peculiar non-minimum phase linear system is presented.
\end{itemize}
Following \cite{mar}, Section \ref{connect} explains the industrial capabilities of classic PIDs by relating them to our intelligent controllers. 
This quite surprising and unexpected result is achieved for the first time to the best of our knowledge.
Section \ref{conclusion} concludes not only by a short list of open problems but also with a discussion of the possible influences on the development of automatic control, which might be 
brought by our model-free standpoint. 

The appendix gives some more explanations on the deduction of Equation \eqref{ultralocal}. 
We are employing
\begin{itemize}
\item rather old-fashioned functional analysis, which goes back to
\cite{volterra1,volterra2,volterra3}. Note that this functional analysis is a mainstay in
engineering since the introduction of {\em Volterra series} (see,
{\it e.g.}, \cite{barrett});
\item some elementary facts stemming from differential algebra
(\cite{kolchin}), which has been quite important in control theory since the
appearance twenty years ago of {\em flatness-based control}
(\cite{flmr}).
\end{itemize}

\section{Model-free control: general principles}\label{mfc}
Our viewpoint on the general principles on model-free control was developed in (\cite{australia,cifa06,ajaccio,esta,malo,milan,marseille}).
\subsection{Intelligent controllers}\label{general}
\subsubsection{Generalities}
Consider again the ultra-local model \eqref{ultralocal} .
Close the loop via the
\emph{intelligent controller}
\begin{equation}\label{ic}
u = - \frac{F - y^{\ast (\nu)} + \mathfrak{C}(e)}{\alpha}
\end{equation}
where
\begin{itemize}
\item $y^\ast$ is the output reference trajectory;
\item $e = y - y^\ast$ is the tracking error;
\item $\mathfrak{C}(e)$ is a \emph{causal}, or \emph{non-anticipative},
functional of $e$, \textit{i.e.}, $\mathfrak{C}(e)$ depends on the past and the present, and not on the future.
\end{itemize}
\begin{remark}
See, \textit{e.g.},
\cite{volterra1,volterra2,volterra3} for an intuitive and clever
presentation of the early stages of the notion of
\emph{functionals}, which were also called sometimes \emph{line
functions}. See Section \ref{app1} in the appendix for more details.
\end{remark}
\begin{remark}
Imposing a reference trajectory $y^\ast$ might lead, as well known, to severe difficulties
with non-minimum phase systems: see, \textit{e.g.}, \cite{predict,sira3,sira} from a flatness-based viewpoint (\cite{flmr,sira}). 
See also Remarks \ref{FLAT}, \ref{remnm}, and Section \ref{conclusion}.
\end{remark}
Combining Equations \eqref{ultralocal} and \eqref{ic} yields the
functional equation
\begin{equation*}\label{fc}
e^{(\nu)} + \mathfrak{C}(e) = 0
\end{equation*}
$\mathfrak{C}$ should be selected such that a perfect tracking is asymptotically ensured,
\textit{i.e.},
\begin{equation}\label{stab}
\lim_{t\rightarrow +\infty} e(t) = 0
\end{equation}
This setting is too general and might not lead to easily
implementable tools. This shortcoming is corrected below.
\subsubsection{Intelligent PIDs}\label{intelpid}
Set $\nu = 2$ in Equation \eqref{ultralocal}. With Equation \eqref{2} define the intelligent
proportional-integral-derivative controller, or iPID, \eqref{ipid}. Combining
Equations \eqref{2} and \eqref{ipid} yields Equation \eqref{stabil}, where $F$ does not appear anymore, {\it i.e.,} the unknown parts and disturbances
of the plant are eliminated. The tracking condition expressed by Equation \eqref{stab} is therefore easily fulfilled by an appropriate tuning of $K_P$, $K_I$, $K_D$. It boils down to 
the stability of a linear differential equation of order $3$, with constant real coefficients.
If $K_I =0$ we obtain an  \emph{intelligent
proportional-derivative controller}, or \emph{iPD},  
\begin{equation}\label{ipd}
\boxed{u = - \frac{F - \ddot{y}^\ast + K_P e + K_D
\dot{e}}{\alpha}}
\end{equation}

Assume now that $\nu = 1$ in Equation \eqref{ultralocal}:
\begin{equation}
\boxed{\dot{y} = F + \alpha u} \label{e1}
\end{equation}
The loop is closed by the \emph{intelligent proportional-integral
controller}, or \emph{iPI},
\begin{equation}\label{ipi}
\boxed{u = - \frac{F - \dot{y}^\ast + K_P e + K_I \int e}{\alpha}}
\end{equation}
Quite often $K_I$ may be set to $0$. It yields an \emph{intelligent
proportional controller}, or \emph{iP},
\begin{equation}\label{ip}
\boxed{u = - \frac{F - \dot{y}^\ast + K_P e}{\alpha}}
\end{equation}
Results in Sections \ref{ord1} and \ref{connect} explain why iPs are
quite often encountered in practice. Their lack of any integration
of the tracking errors demonstrate that the anti-windup algorithms,
which are familiar with ``classic'' PIDs and PIs, are no more
necessary.
\begin{remark}
There is, as
well known, a huge literature on ``classic'' PIDs and PIs in order to give efficient rules for the gain tuning. Those recipes are too often rather intricate. See,
\textit{e.g.}, the two books by \cite{astrom1}, \cite{od}, and the numerous references therein.
\end{remark}
\begin{remark}\label{FLAT}
Output reference trajectories of the form $y^\ast$
do not seem to be familiar in industrial applications of classic PIDs. This
absence often leads to disturbing oscillations, and mismatches like overshoots and undershoots. Selecting $y^\ast$ plays of course a key r\^ole in the
implementation of the control synthesis. Mimicking for this tracking
the highly effective feed-forward flatness-based viewpoint (see,
{\it e.g.}, \cite{flmr}, and \cite{levine,sira}, and the numerous
references in those two books) is achieved in Section
\ref{restrict} where a part of the system, which happens to be flat,
is already known. This is unfortunately impossible in general: are
systems like \eqref{functional} and/or \eqref{eq} in the appendix flat or not? Even
if the above systems were flat, it might be difficult then to verify
if $y$ is a flat output or not. 
\end{remark}
\begin{remark}
For obtaining a
suitable trajectory planning, impose to $y$ to satisfy a given ordinary
differential equation. It permits moreover if
the planning turns out to be poor because of some abrupt change to
replace quite easily the preceding equation by another one.
\end{remark}

\subsection{Other possible intelligent controllers}\label{extension}
The \emph{generalized proportional-integral} controllers, or \emph{GPIs}, were introduced by \cite{GPI} in order to tackle some tricky problems like those stemming 
from non-minimum phase systems. Several practical case-studies have confirmed their usefulness (see, \textit{e.g.}, \cite{siragpi,moragpi}). Although it would be possible to define their 
\emph{intelligent} counterparts in general, we are limiting ourselves here to a single case which will be utilized in Section \ref{nm}. Replace the ultra-local  model \eqref{e1} by
\begin{equation}\label{2i}
\dot y=F+\alpha u+\beta \int u
\end{equation}  
where $\alpha, \beta \in \mathbb{R}$ are constant. Set in Equation \eqref{ic}
\begin{equation}\label{i2}
\mathfrak{C}(e) =K_Pe+K_I\int e+K_{II}\int\int e
\end{equation}
where $K_I, K_{II} \in \mathbb{R}$ are suitable constant gains. See Section \ref{2int} for an analogous regulator.

\section{Online estimation of $F$}\label{2}
Our first publications on model-free
control were proposing for the estimation of $F$ recent techniques on the
numerical differentiations of noisy signals (see \cite{nl}, and \cite{mboup,liu}) for estimating $y^{({\nu})}$ in Equation \eqref{ultralocal}. Existing applications were until today based
on a simple version of this differentiation procedure, which is quite close to what is presented in this Section, 
namely the utilization of the parameter identification techniques by \cite{sira1,sira2}. 

\subsection{General principles}
The approximation of an integrable function, \textit{i.e.}, of a quite general function $[a, b] \rightarrow  \mathbb{R}$, $a, b \in \mathbb{R}$, $a < b$, by a \emph{step} function $F_{\text{approx}}$, \textit{i.e.}, a piecewise constant 
function, is classic in mathematical analysis (see, {\it e.g.}, the excellent textbooks by \cite{godement} and \cite{rudin1}). A suitable approximate estimation of $F$ in Equation \eqref{ultralocal} boils down therefore to the estimation 
of the constant parameter $\phi$ in Equation \eqref{identif} if it can be achieved during a sufficiently ``small'' time interval. Analogous estimations of $F$ may be carried on via the intelligent 
controllers \eqref{ipid}-\eqref{ipd}-\eqref{ipi}-\eqref{ip}.

\subsection{Identifiability via operational calculus}\label{review}
\subsubsection{Operatational calculus}
In order to encompass all the previous equations, where $F$ is replaced by $F_{\text{approx}}$, consider the equation, where 
the classic rules of operational calculus are utilized (\cite{miku,yosida}),
\begin{equation}
L_1(s)Z_1 + L_2(s)Z_2 = \frac{\phi}{s} +  I(s) \label{II}
\end{equation}
\begin{itemize}
\item $\phi$ is a constant real parameter, which has to be identified;
\item $L_1, L_2 \in \mathbb{R}[s, s^{-1}]$ are Laurent polynomials;
\item $I \in \mathbb{R}[s]$ is a polynomial associated to the initial conditions.
\end{itemize}
Multiplying both sides of Equation \eqref{II} by $\frac{d^N
}{ds^N}$, where $N$ is large enough, permits to get rid of the
initial conditions. It yields the \emph{linear identifiability}
(\cite{sira1,sira2}) of $\phi$ thanks to the formula
\begin{equation}\label{III}
\frac{(-1)^N N!}{s^{N+1}} \phi = \frac{d^N }{ds^N} \left( L_1(s)Z_1
+ L_2(s)Z_2 \right)
\end{equation}
Multiplying both sides of Equation \eqref{III} by $s^{-M}$, where $M
> 0$ is large enough, permits to get rid of positive powers of
$s$, \textit{i.e.}, of derivatives with respect to time. 
\begin{remark}
Sometimes it might be interesting in practice to replace $s^{-M}$ by a suitable rational function of $s$, \textit{i.e.}, by a suitable element of  $\mathbb{R}(s)$.
\end{remark}
\subsubsection{Time domain}\label{td}
The remaining negative
powers of $s$ correspond to iterated time integrals. The corresponding formulae in
the time domain are easily deduced thanks to the correspondence
between $\frac{d^\kappa}{ds^\kappa}$, $\kappa \geq 1$, and the multiplication by $(-
t)^\kappa$ in the time domain (see some examples in Section \ref{isch}). They may be easily implemented as
discrete linear filters.

\subsection{Noise attenuation}
The notion of noise, which is usually described in engineering and, more generally, in applied sciences via probabilistic and 
statistical tools, is borrowed here from \cite{bruit} (see also \cite{lobry}, and the references therein on \emph{nonstandard analysis}). 
Then the noise is related to \emph{quick fluctuations} around zero.  Such a fluctuation
is a Lebesgue-integrable real-valued time function $\mathcal{F}$ which is characterized by the following property: 

\noindent its integral $\int_{\tau_i}^{\tau_f} \mathcal{F}(\tau) d\tau$ over any finite
interval is \emph{infinitesimal}, \textit{i.e.}, very ``small''. 

\noindent The robustness with respect to corrupting noises is thus explained thanks to Section \ref{td}. 

\begin{remark}
This standpoint on denoising has not only been confirmed by several applications of model-free control, which were already cited in the introduction, but 
also by numerous ones in model-based linear control and in signal processing (see, \textit{e.g.}, \cite{rupture,at,morales,pereira1,pereira2,trapero1,trapero2,trapero3}).
Note moreover that the nonlinear estimation techniques advocated by \cite{nl} exhibit for the same reason ``good'' robustness properties, which were already illustrated by several 
case-studies (see, \textit{e.g.}, \cite{menhour,mora}, and the references therein). 
\end{remark}

\subsection{Some more explicit calculations}\label{isch}
\subsubsection{First example}
With Equation \eqref{e1}, Equation \eqref{II} becomes
$$ s Y = \frac{\phi}{s} + \alpha U + y_0 $$
where
\begin{itemize}
\item $y_0$ is the initial condition corresponding to the time
interval $[t - L, t]$,
\item $\phi$ is a constant.
\end{itemize}
Get rid of $y_0$ by multiplying both sides by $\frac{d}{ds}$:
$$y+s\frac{dy}{ds}=-\frac{\phi}{s^2}+\alpha \frac{du}{ds}$$
Multiplying both sides by $s^{-2}$ for smoothing the noise yields in
time domain yields
$$\phi = -\frac{6}{L^3}\int_{t-L}^t
\left((L-2\sigma)y(\sigma)+\alpha\sigma(L-\sigma)u(\sigma)
\right)d\sigma$$
where $L$ is quite small.
\begin{remark}
$L$ depends of course on
\begin{itemize}
\item the sampling period,
\item the noise intensity.
\end{itemize}
Both may differ a lot as demonstrated by the numerous references on concrete case-studies given at the beginning of the introduction.
\end{remark}

\subsubsection{Second example}
Close the loop with the iP \eqref{ip}. It yields
$$
\phi = \frac{1}{L}\left[\int_{t-L}^{t}\left(\dot{y}^{\star}-\alpha u
- K_P e \right) d\sigma \right] 
$$

\section{When is the order $\nu = 1$ enough?}\label{ord1} 
\begin{figure}[htp]
\centering
\includegraphics[width=9.05cm]{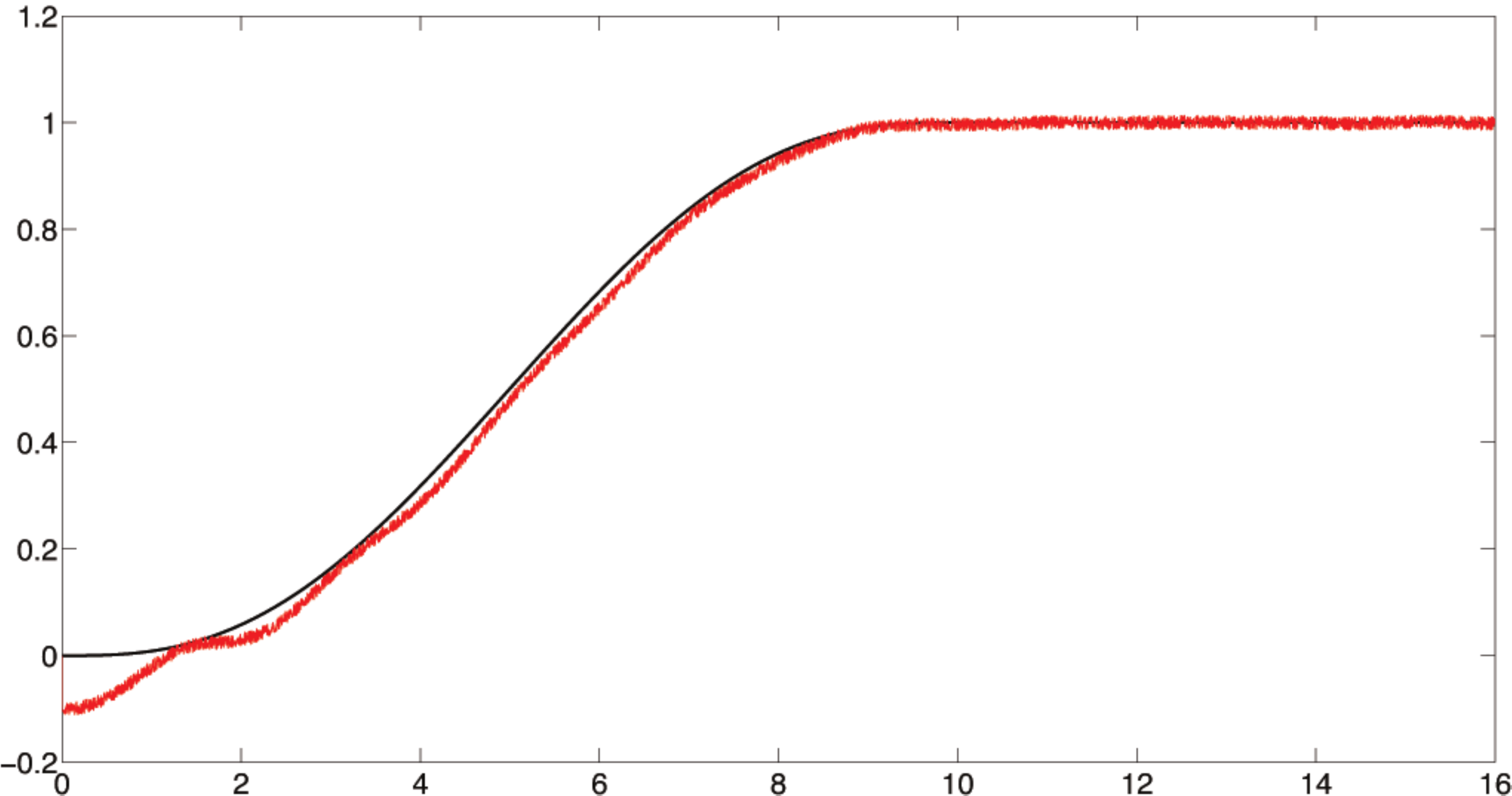}
\caption{System output and reference}\label{fig01}
\end{figure}

\begin{figure}[htp]
\centering
\includegraphics[width=9.05cm]{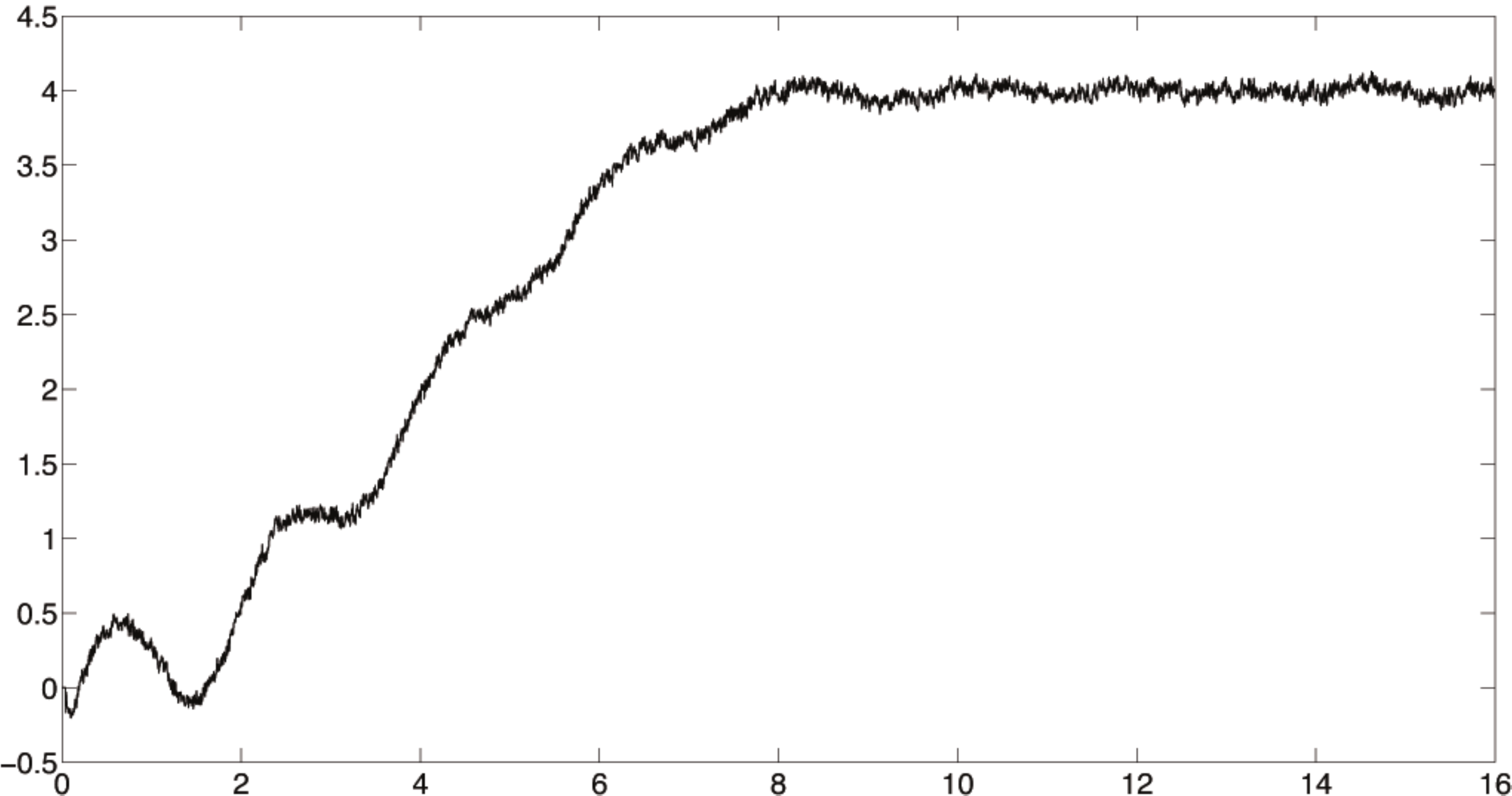}
\caption{iPI control}\label{fig02}
\end{figure}

\begin{figure}[htp]
\centering
\includegraphics[width=9.05cm]{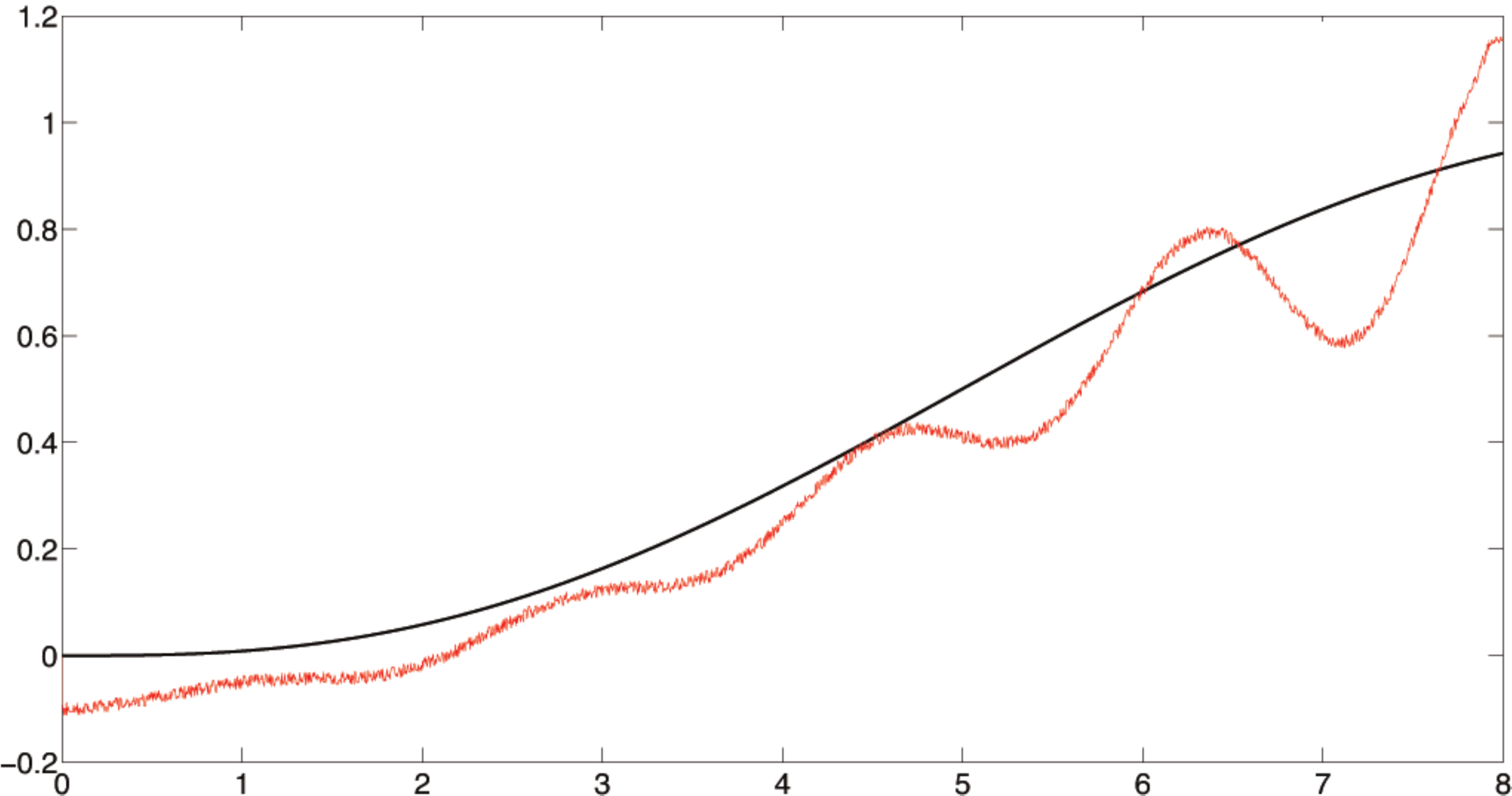}
\caption{System output and reference}\label{fig03}
\end{figure}

A most
notable exception in the choice of a first order ultra-local model,
\textit{i.e.}, $\nu = 1$ in Equation \eqref{ultralocal}, is provided
by the magnetic bearing studied by \cite{cifa12}, where the friction
is almost negligeable. Start therefore with the elementary constant linear
system
\begin{equation}
\ddot{y}+c\dot{y}+4 y=u \label{oscillateur}
\end{equation}
where $c\dot{y}$ stands for some elementary friction. Figures
\ref{fig01} and \ref{fig02} yield satisfactory numerical simulations
with a iPI corrector. The following values were selected for the
parameters: $c=3$, $\alpha = 1$, $K_P=16$, $K_I=25$. With a harmonic
oscillator, where $c = 0$, Figure \ref{fig03} displays on the other
hand a strong degradation of the performances with an iPI. Lack of
friction in a given system might be related to the absence of
$\dot{y}$ in the unknown equation. Taking $\nu = 1$ in Equation
\eqref{ultralocal} would therefore yield an ``algebraic loop,''
which adds numerical instabilities and therefore deteriorates the
control behavior.


\section{Numerical experiments}\label{ex}
In the subsequent simulations the sampling time is $T_e=0.01 \text{s}$. The
corrupting noise is additive, normal, zero-mean, with a standard
deviation equal in Sections \ref{restrict} and \ref{fourier} to $0.01$, and to $0.03$ elsewhere.

\subsection{Control with a partially known system}\label{restrict}

\subsubsection{A crude description}\label{crude}

\begin{figure}
\begin{center}
\subfigure[Tustin's model]{\resizebox*{5cm}{!}{\includegraphics{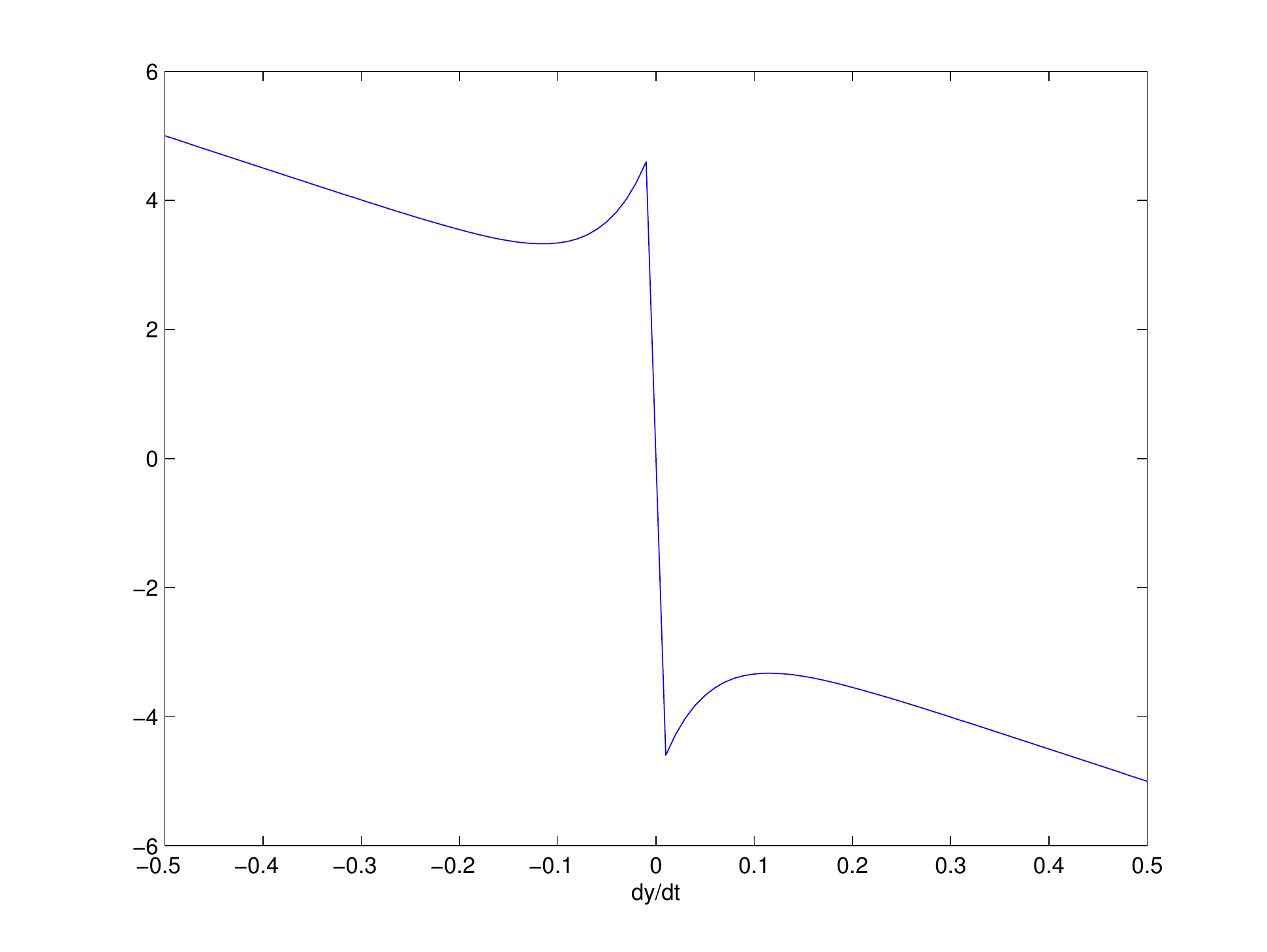}}}%
\subfigure[Friction with the iPID]{\resizebox*{5cm}{!}{\includegraphics{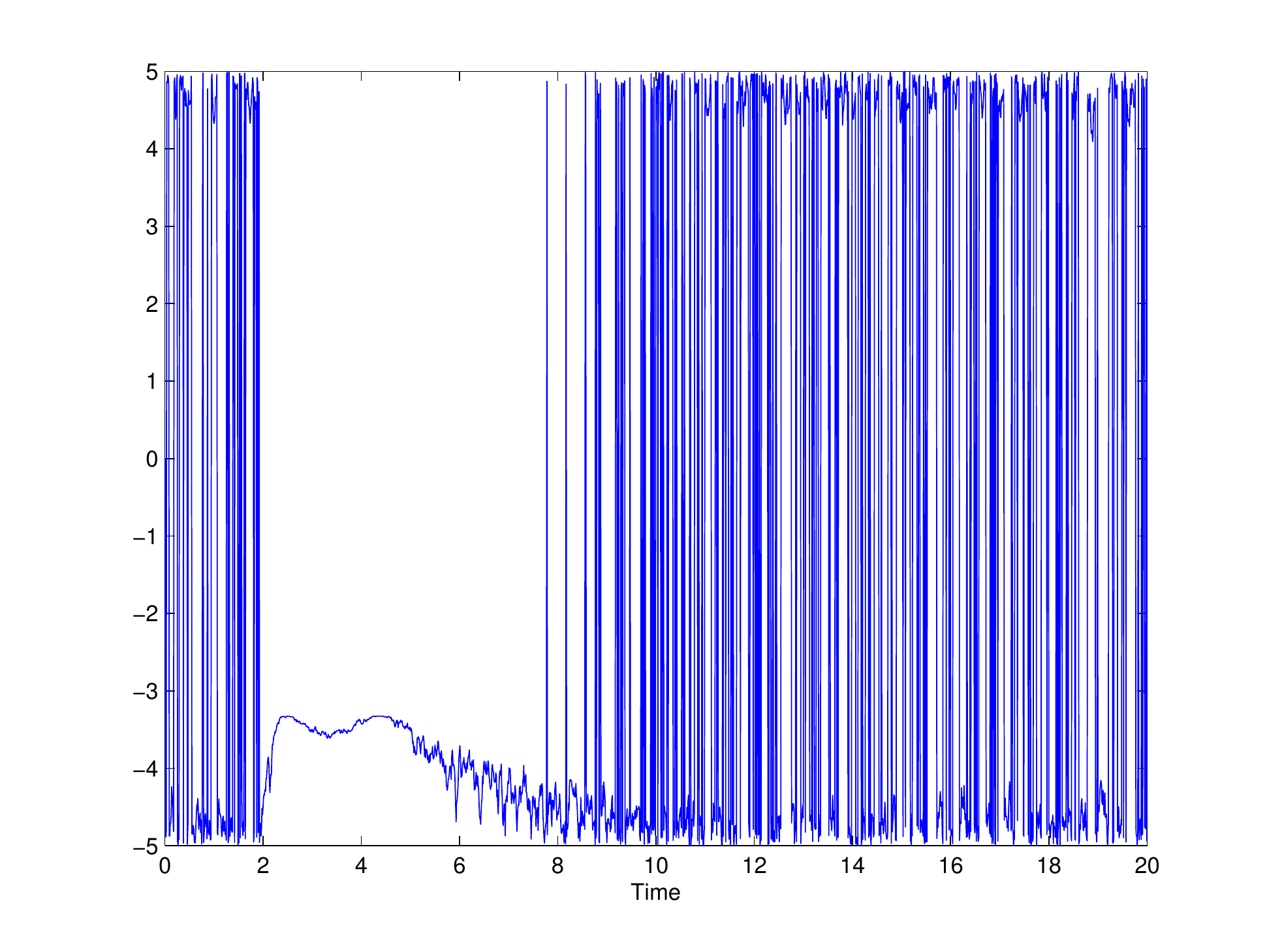}}}%
\caption{Model and time evolution of friction}%
\label{6forces}
\end{center}
\end{figure}

Consider the non-linear Duffing spring with friction:
\begin{equation}\label{ressort}m\ddot{y}= -{\mathcal{K}}(y) +
{\mathcal{F}}(\dot{y})-d\dot y+F_{\tiny{\mbox{\rm
ext}}}\end{equation} where
\begin{itemize}
\item $y$ is the length of the spring,
\item $m$ is a point mass,
\item $F_{\tiny{\mbox{\rm ext}}}=u$ is the control variable,
\item ${\mathcal{K}}(y) = k_1 y + k_3 y^3$ is the resulting force from the Hooke law and the Duffing cubic term,
\item $d {\dot y}$ is a classic linear friction and ${\mathcal{F}}(\dot y)$ a non-linear one.  The term
${\mathcal{F}}(\dot y)$ corresponds to  the Tustin friction
(\cite{tustin}) (see, also, \cite{ols}), which is rather
violent with respect to the sign change of the speed (see Figure
\ref{6forces}).
\end{itemize}
Set $m=0.5$,  $k_1 = 3$. The partially known system
$$
m\ddot y= k_1 y + u
$$
is flat, and $y$ is a flat output. It helps us to determine a
suitable reference trajectory $y^\star$ and the corresponding
nominal control variable $u^\star = m\ddot y^\star + k_1 y^\star$.
In the numerical simulations, we utilize $\hat{k}_1 = 2$, $d = 1$,
$k_3 = 2$, which are in fact unknown.

\subsubsection{A PID controller}\label{classpid}
Set $u = u^\star  + v$. Associate to $v$ a PID corrector for
alleviating the tracking error $e=y-y^\star$ by imposing a
denominator of the form $(s + 1.5)^3$.  The corresponding tuning
gains are $k_P = 1.375$, $k_I =1.6875$, $k_D = 2.25$.

\subsubsection{iPID}
The main difference of the iPID  is the following one: The presence
of $F$ which is estimated in order to compensate the nonlinearities
and the perturbations like frictions. For comparison purposes, its
gains are the same as previously.

\subsubsection{iP}
We do not take any advantage of Equation \eqref{ressort}. The error
tracking dynamics is again given with a pole equal to $- 1.5$, {\it i.e.}, by the denominator $(s+1.5)$.
\subsubsection{Numerical experiments}

\begin{figure}
\begin{center}
\subfigure[Controls $F_{{\rm ext}}(t)$: PID(--, blue), iP(--, red), iPID(--, green)]{
\resizebox*{15cm}{!}{\includegraphics{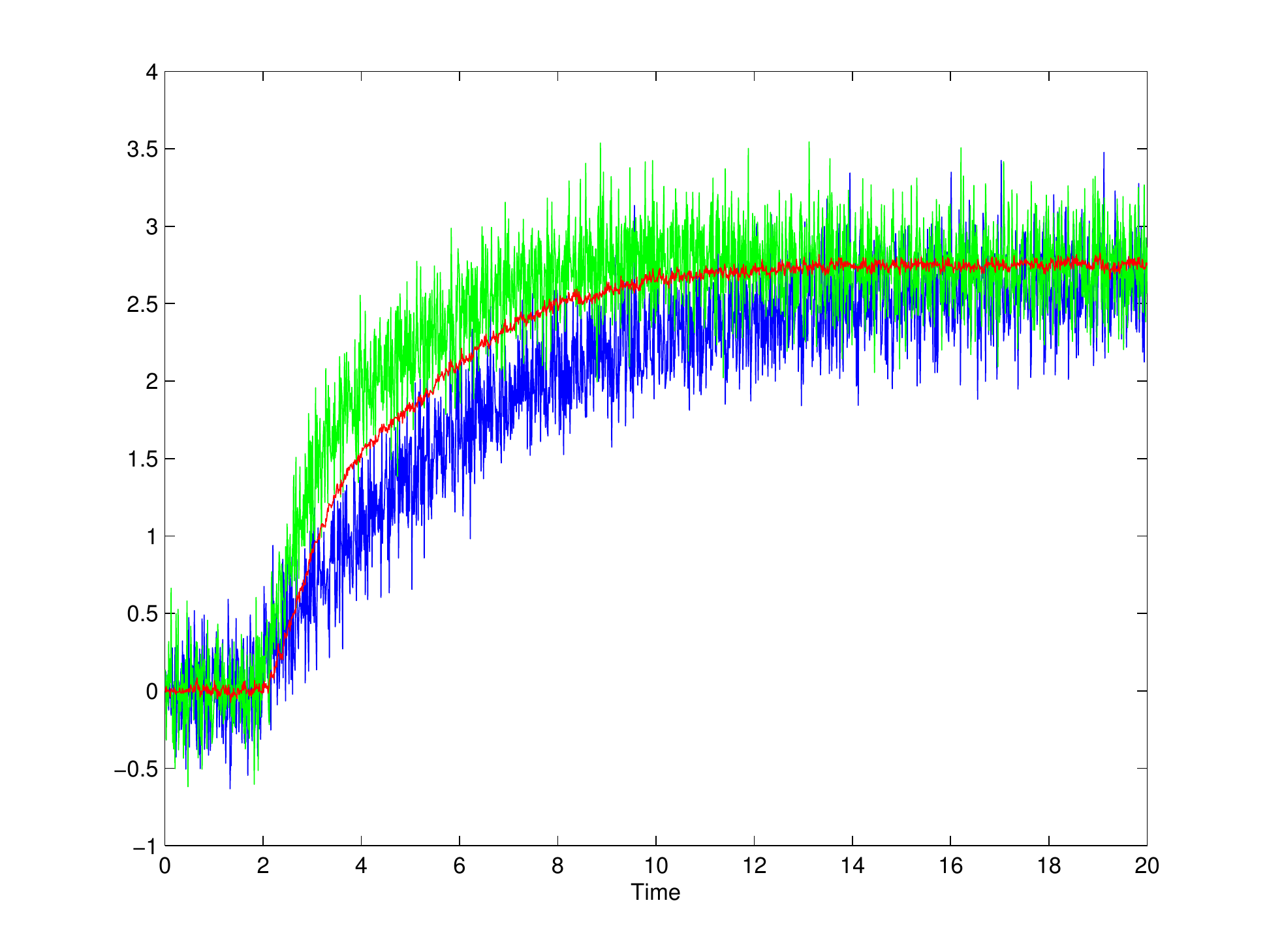}}}\\%
\subfigure[Setpoint (- ., black), reference (- -, black), and output: PID(--, blue), iP(--, red), iPID(--, green)]{
\resizebox*{15cm}{!}{\includegraphics{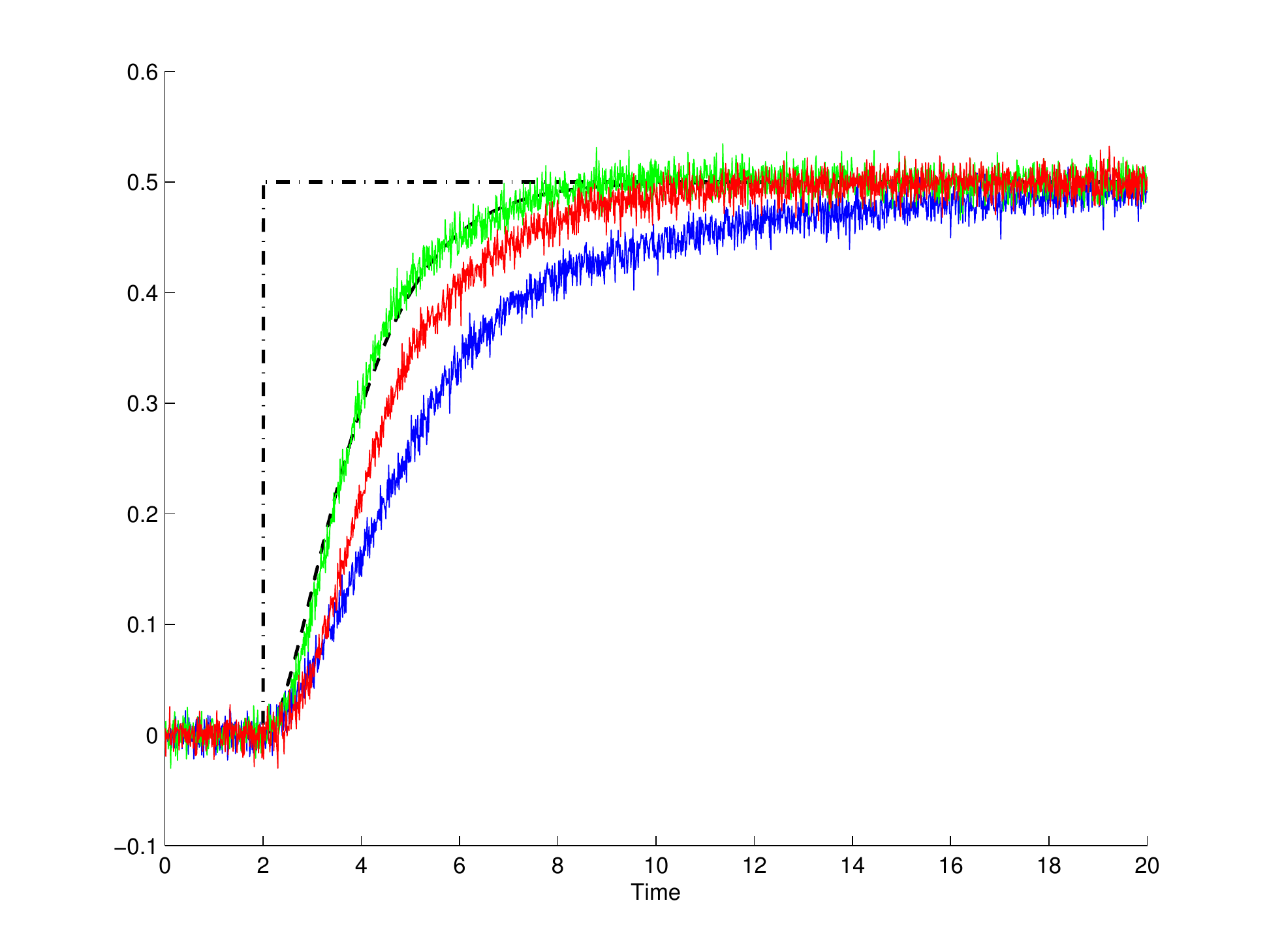}}}%
\caption{Partially known system: comparisons}%
\label{6tout}
\end{center}
\end{figure}

Figure \ref{6tout} shows quite poor results with the PID of Section
\ref{classpid} . They become excellent with the iPID, and correct
with iP. The practician might be right to prefer this last control
synthesis
\begin{itemize}
\item where the implementation is immediate,
\item if a most acute
precision may be neglected.
\end{itemize}

\begin{remark}
See also \cite{vil1} and \cite{vil2} for concrete examples related to guided vehicles.
\end{remark}

\subsection{Robustness with respect to system's changes}
The examples below demonstrate that if the system is changing, our
intelligent controllers behave quite well without the need of any
new calibration.

\subsubsection{Scenario 1: the nominal case}
\begin{figure}
\begin{center}
\subfigure[Controls: PID(--, blue), iP(--, red)]{
\resizebox*{15cm}{!}{\includegraphics{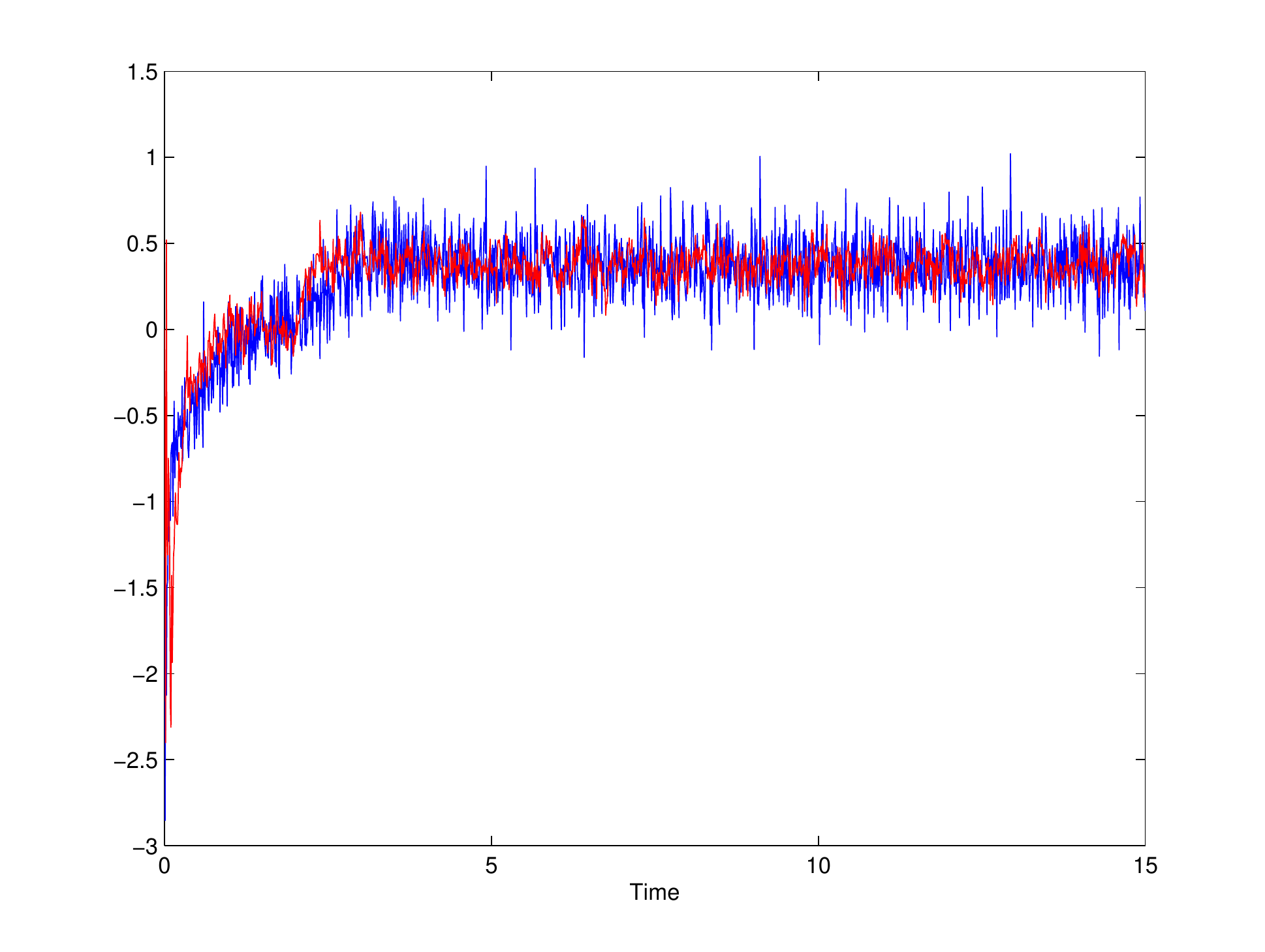}}}\\
\subfigure[Setpoint (- ., black), reference (- -, black), and outputs:  PID(--, blue), iP(--, red)]{
\resizebox*{15cm}{!}{\includegraphics{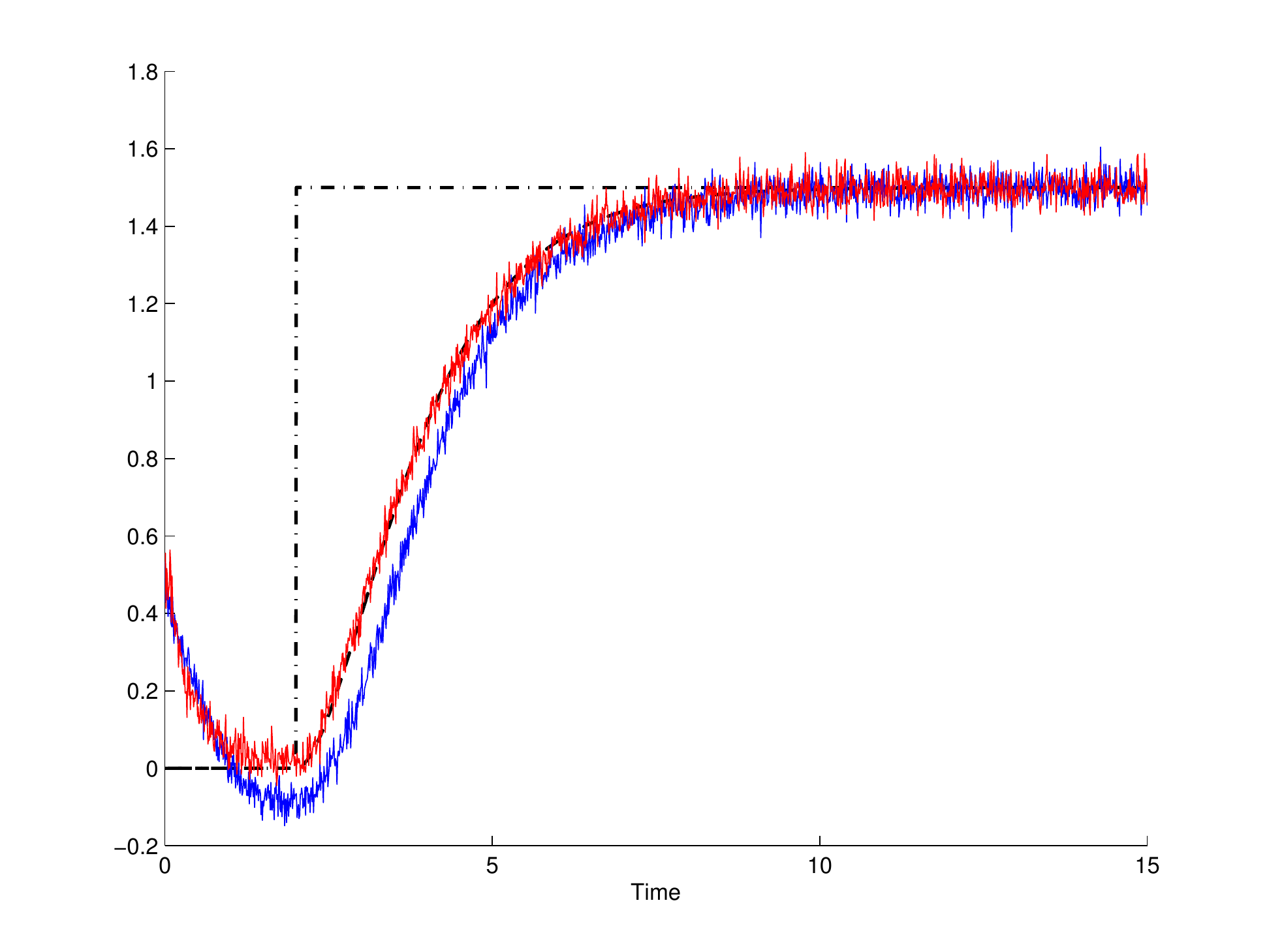}}}%
\caption{Scenario 1: comparisons}%
\label{1XS1}
\end{center}
\end{figure}

The nominal system is defined by the transfer function
\begin{equation}\label{transfer}
\frac{(s+2)^2}{(s+1)^3}
\end{equation}
A tuning of a classic PID
controller
\begin{equation}\label{classic}
u = k_p e + k_i \int e + k_d \dot{e}
\end{equation}
where
\begin{itemize}
\item  $e = y - y^\ast$ is the tracking error,
\item $k_p, k_i, k_d \in \mathbb{R}$ are the gains,
\end{itemize}
yields via standard techniques (see, \textit{e.g.}, \cite{astrom1})
$k_p=1.8177$, $k_i=0.7755$, $k_d=0.1766$. A low-pass filter is
moreover added to the derivative $\dot{e}$ in order to attenuate the
corrupting noise. Our model-free approach utilizes the ultra-local
model $\dot{y} = F + u$ and an iP \eqref{ip} where $K_P= 1.8177$.
Figure \ref{1XS1} shows perhaps a slightly better
behavior of the iP.

\subsubsection{Scenario 2: modifying the pole}\label{2.2}
\begin{figure}
\begin{center}
\subfigure[Controls: PID(--, blue), iP(--, red)]{
\resizebox*{15cm}{!}{\includegraphics{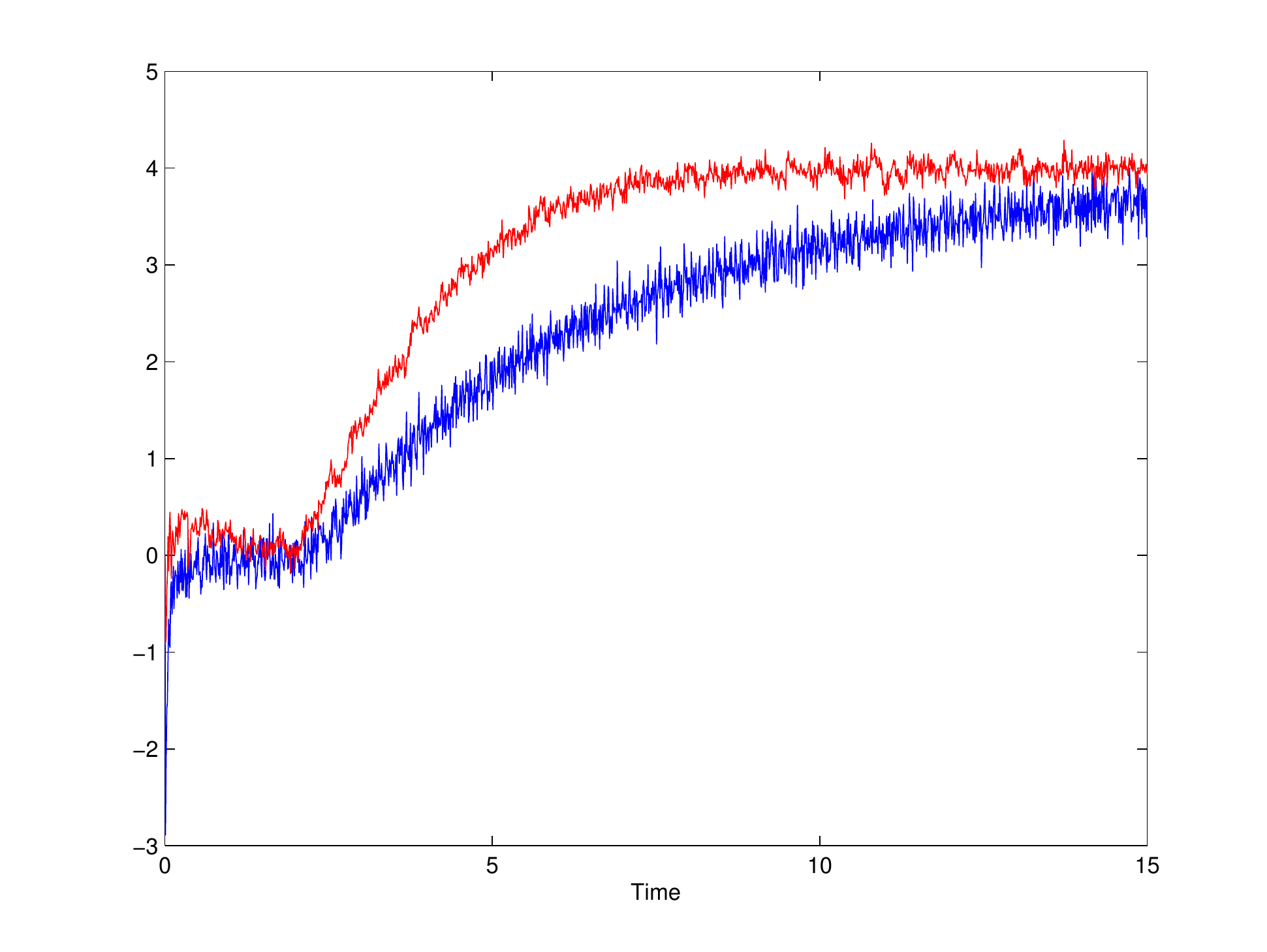}}}\\
\subfigure[Setpoint (- ., black), reference (- -, black), and outputs:  PID(--, blue), iP(--, red)]{
\resizebox*{15cm}{!}{\includegraphics{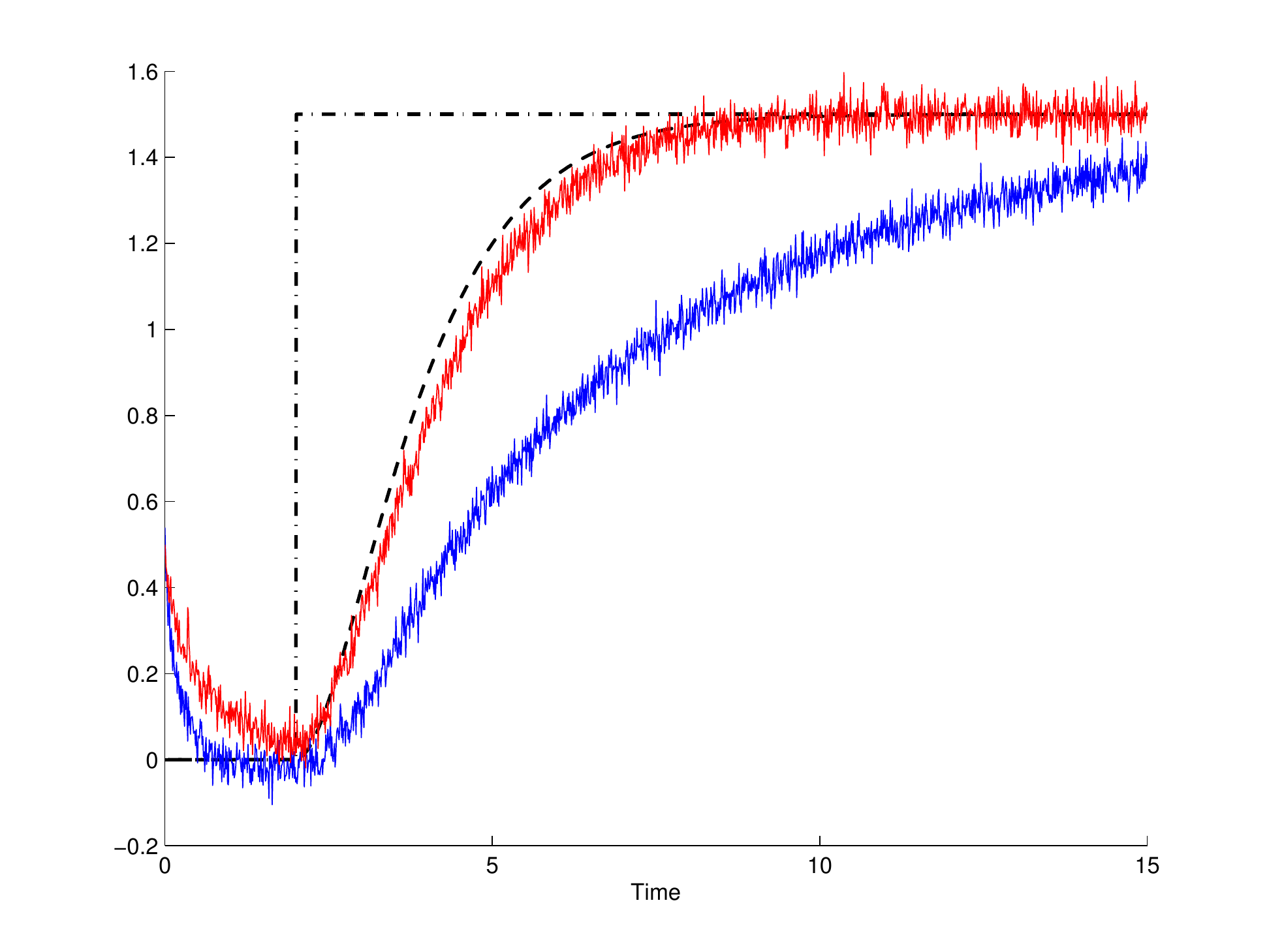}}}%
\caption{Scenario 2: comparisons}%
\label{1XS2}
\end{center}
\end{figure}

A system change, aging for instance, might be seen by as new pole
$-2.2$ in the transfer function \eqref{transfer} which becomes
\begin{equation*}\label{transfer2}
\frac{(s+2)^2}{(s+2.2)^3}
\end{equation*}
As shown in Figure \ref{1XS2}, without any new calibration the performances of the PID worsen whereas
those of the iP remain excellent.

%

\subsubsection{Scenario 3: actuator's fault}\label{8}
\begin{figure}
\begin{center}
\subfigure[Controls: PID(--, blue), iP(--, red)]{
\resizebox*{15cm}{!}{\includegraphics{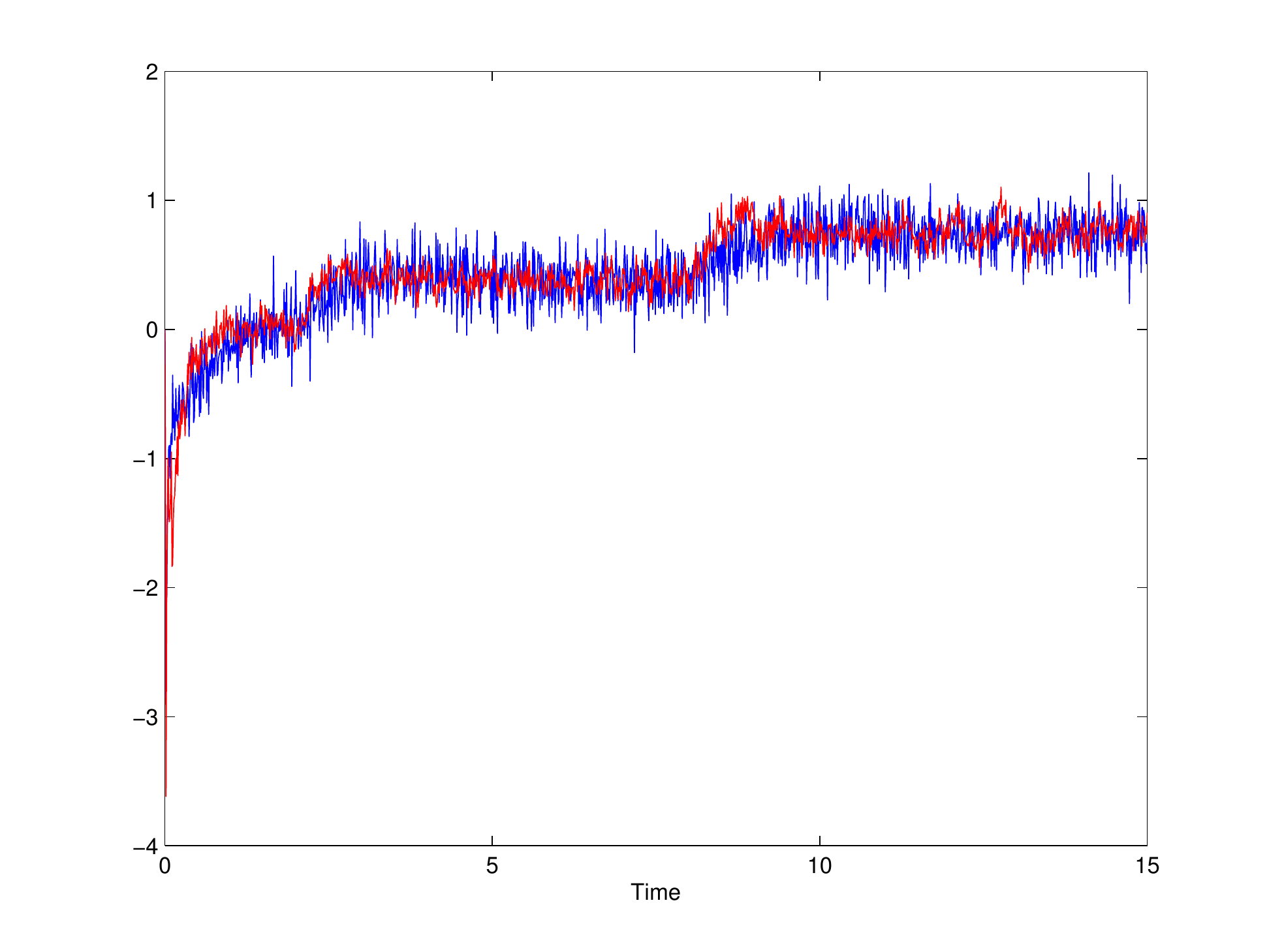}}}\\
\subfigure[Setpoint (- ., black), reference (- -, black), and outputs:  PID(--, blue), iP(--, red)]{
\resizebox*{15cm}{!}{\includegraphics{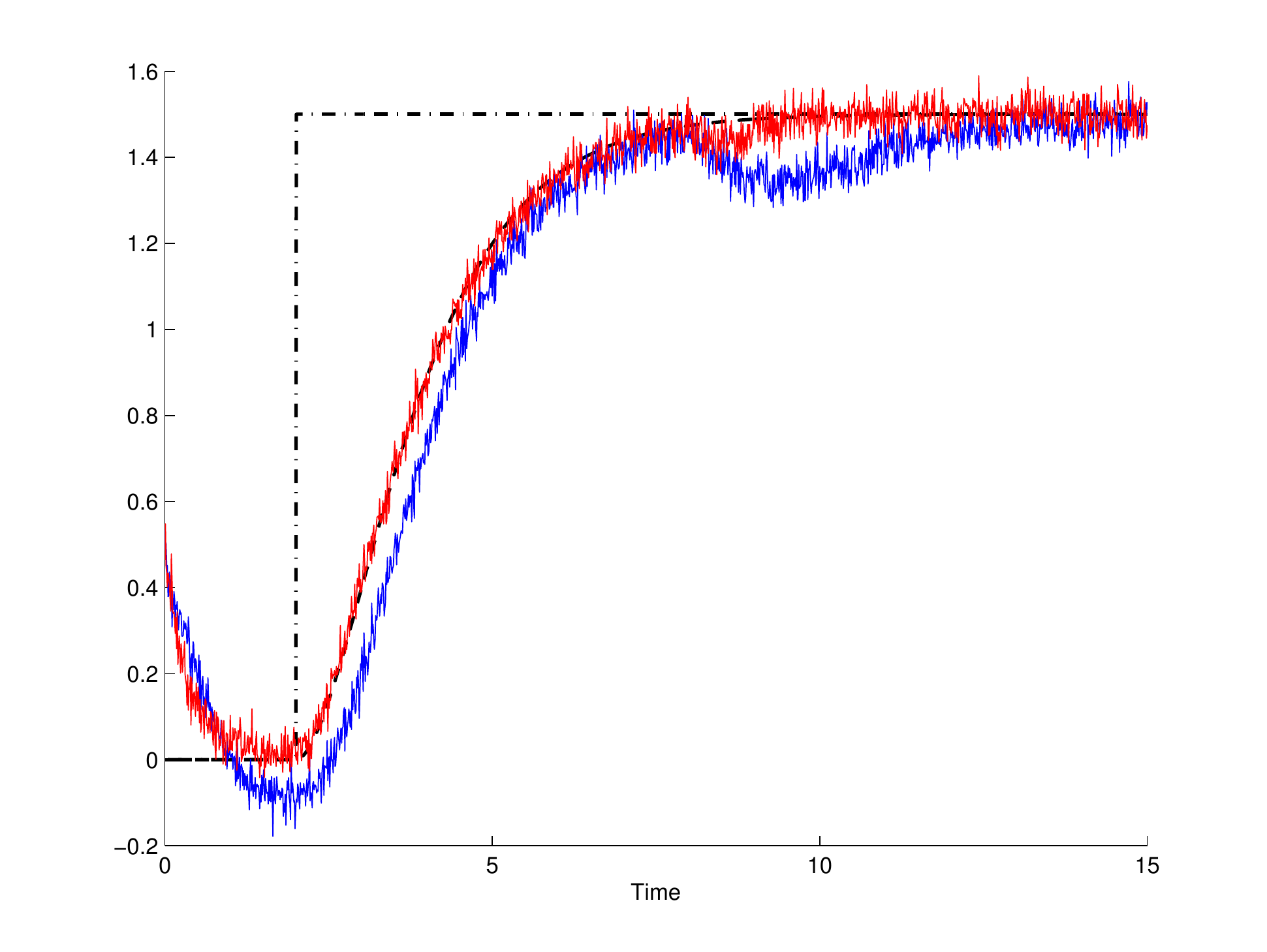}}}%
\caption{Scenario 3: comparisons}%
\label{1XS3}
\end{center}
\end{figure}

A power loss of the actuator occurs at time $t=8\text{s}$. It is
simulated by dividing the control by $2$ at $t=8\text{s}$. Figure
\ref{1XS3} shows an accommodation of the iP which is
much faster than with the PID.

\begin{remark}
Sections \ref{2.2} and \ref{8} may be understood as instances of \emph{fault accommodation}, which contrarily to most
of the existing literature are not model-based (see, also,
\cite{moussa}). It is perhaps worth mentioning here that
model-based fault diagnosis has also benefited from the estimation
techniques summarized in Section \ref{2} (see \cite{diag,nl}).
\end{remark}

\subsection{A non-linear system}
\begin{figure}
\begin{center}
\subfigure[Controls: PID(--, blue), iP(--, red)]{
\resizebox*{15cm}{!}{\includegraphics{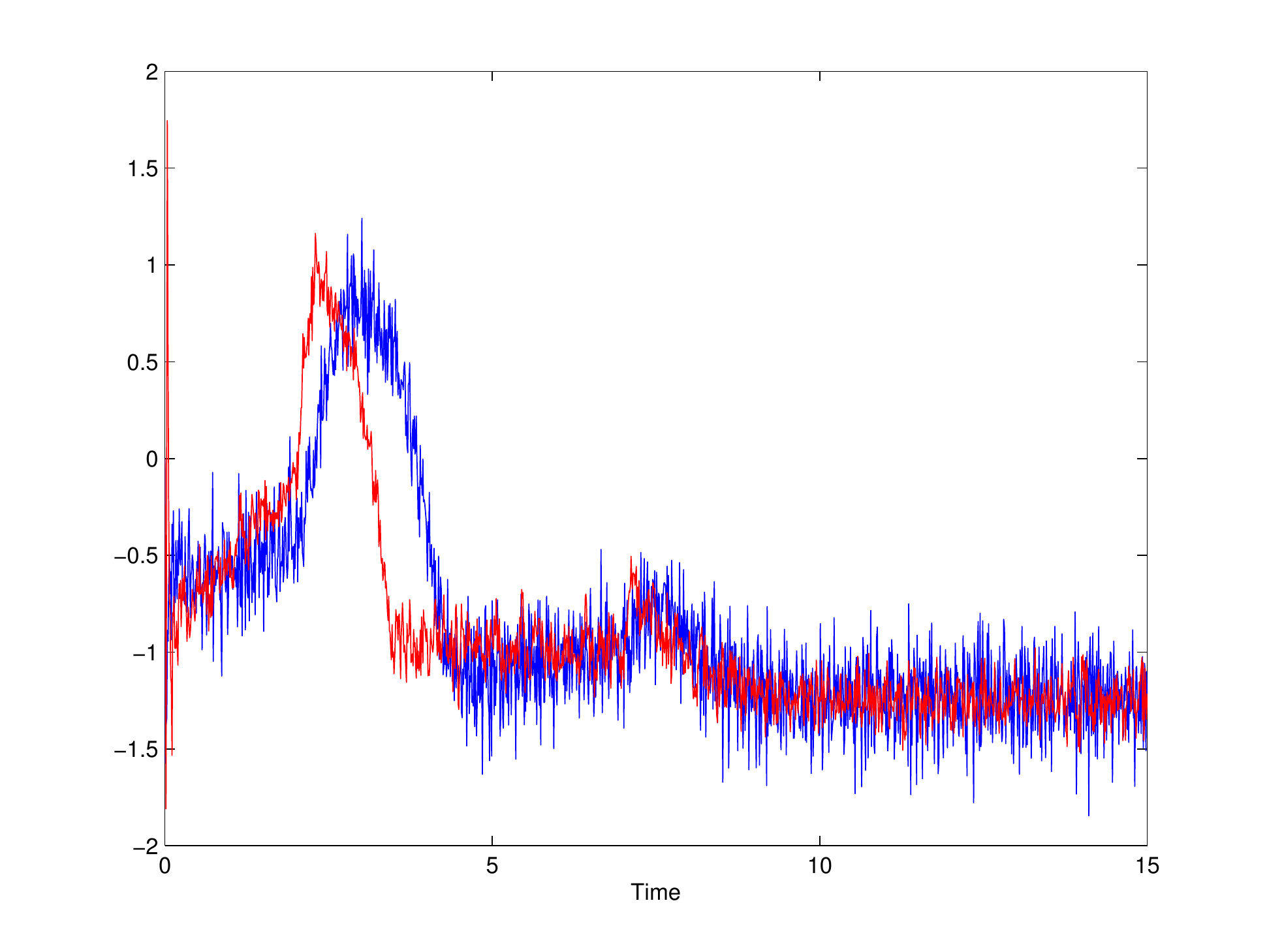}}}\\
\subfigure[Setpoint (- ., black), reference (- -, black), and outputs:  PID(--, blue), iP(--, red)]{
\resizebox*{15cm}{!}{\includegraphics{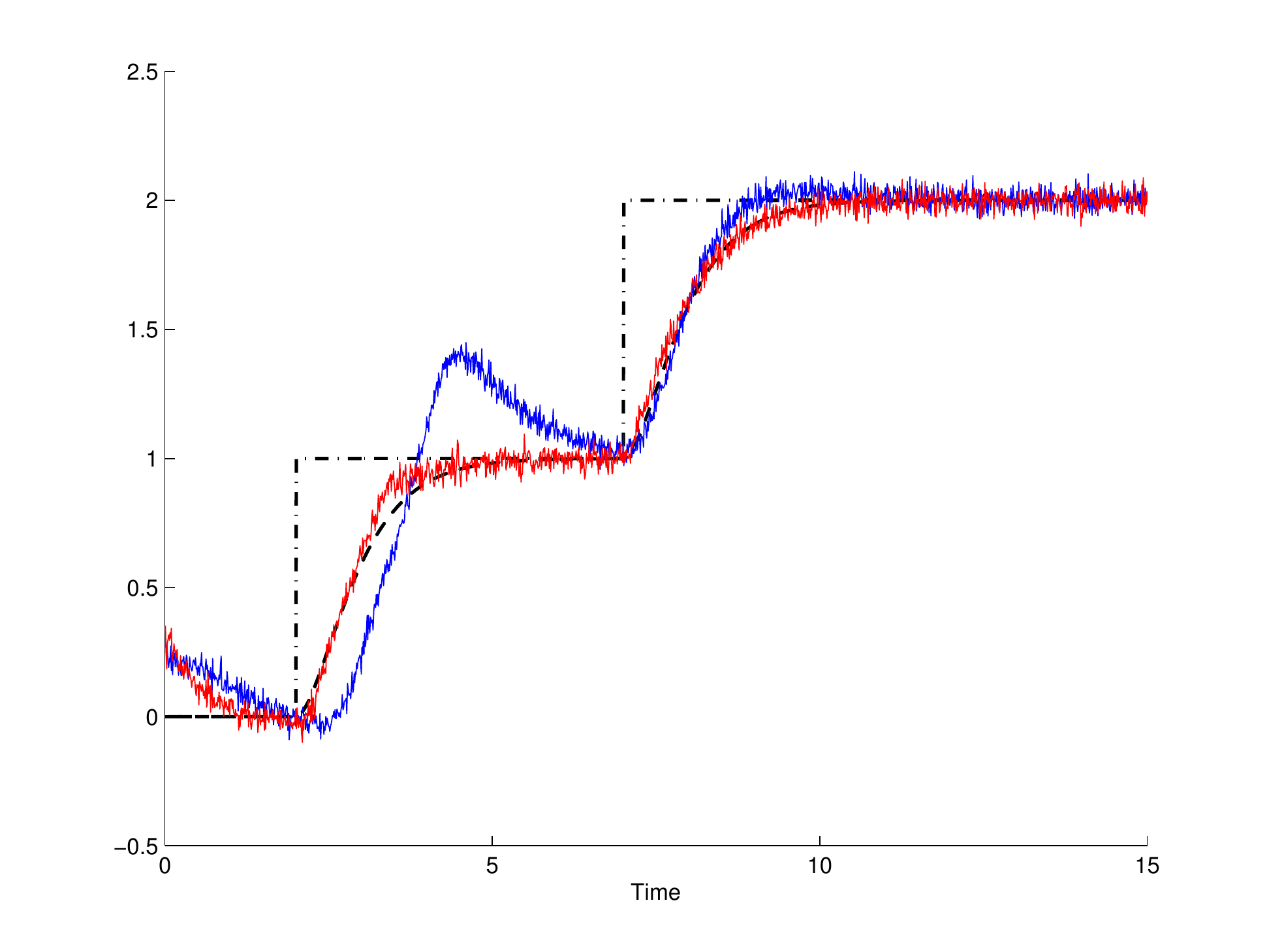}}}%
\caption{Non-linear system: comparisons}%
\label{2X}
\end{center}
\end{figure}

Take the following academic unstable non-linear system
$$\dot y-y=u^3$$
The clssic PID \eqref{classic} is tuned with $k_p=2.2727$,
$k_i=1.8769$, $k_d=0.1750$. The simulations depicted in Figure
\ref{2X} shows a poor trajectory tracking for small values of the
reference trajectory. The iP, which is related to the ultra-local
model $\dot y=F+u$, corresponds to $K_P=2.2727$. Its excellent
performances in the whole operating domain are also shown in Figure
\ref{2X}.

%

\subsection{Delay systems}\label{del}
\begin{figure}
\begin{center}
{\resizebox*{5cm}{!}{\includegraphics{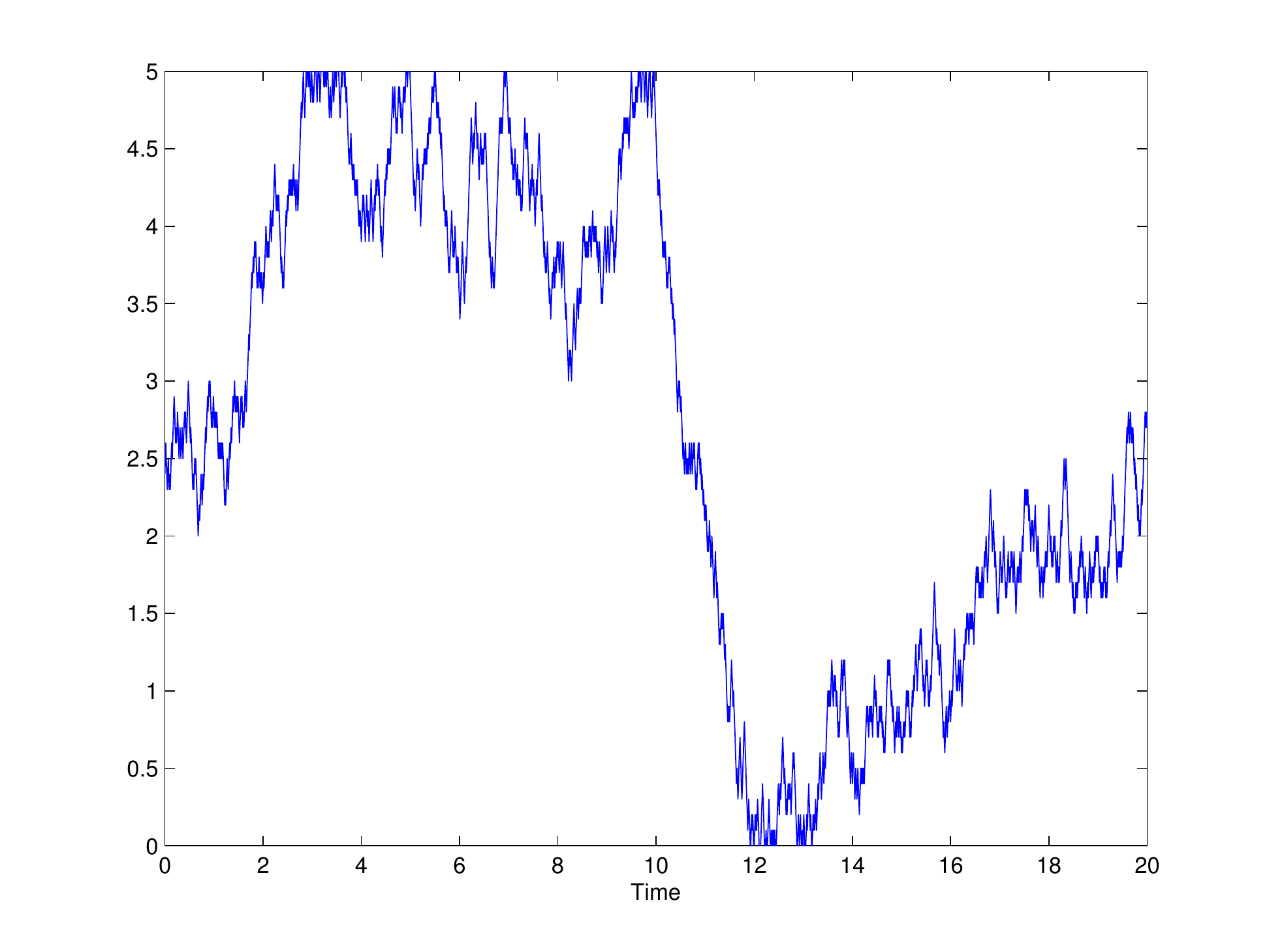}}}%
\caption{Delay function}%
\label{3XD}
\end{center}
\end{figure}

\begin{figure}
\begin{center}
\subfigure[Control]{
\resizebox*{15cm}{!}{\includegraphics{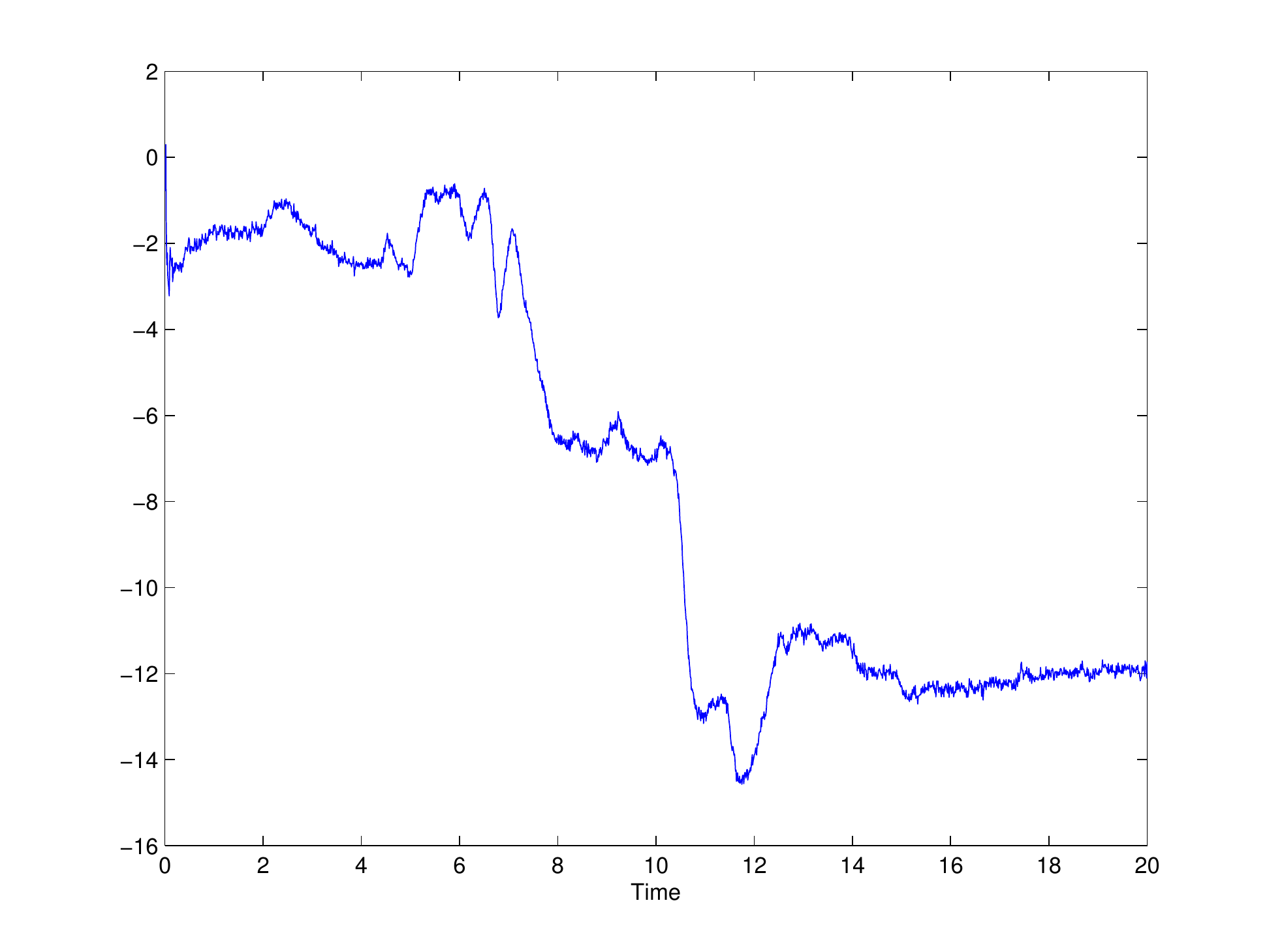}}}\\
\subfigure[Setpoint (- ., black), reference (- -, black), and output (--, blue)]{
\resizebox*{15cm}{!}{\includegraphics{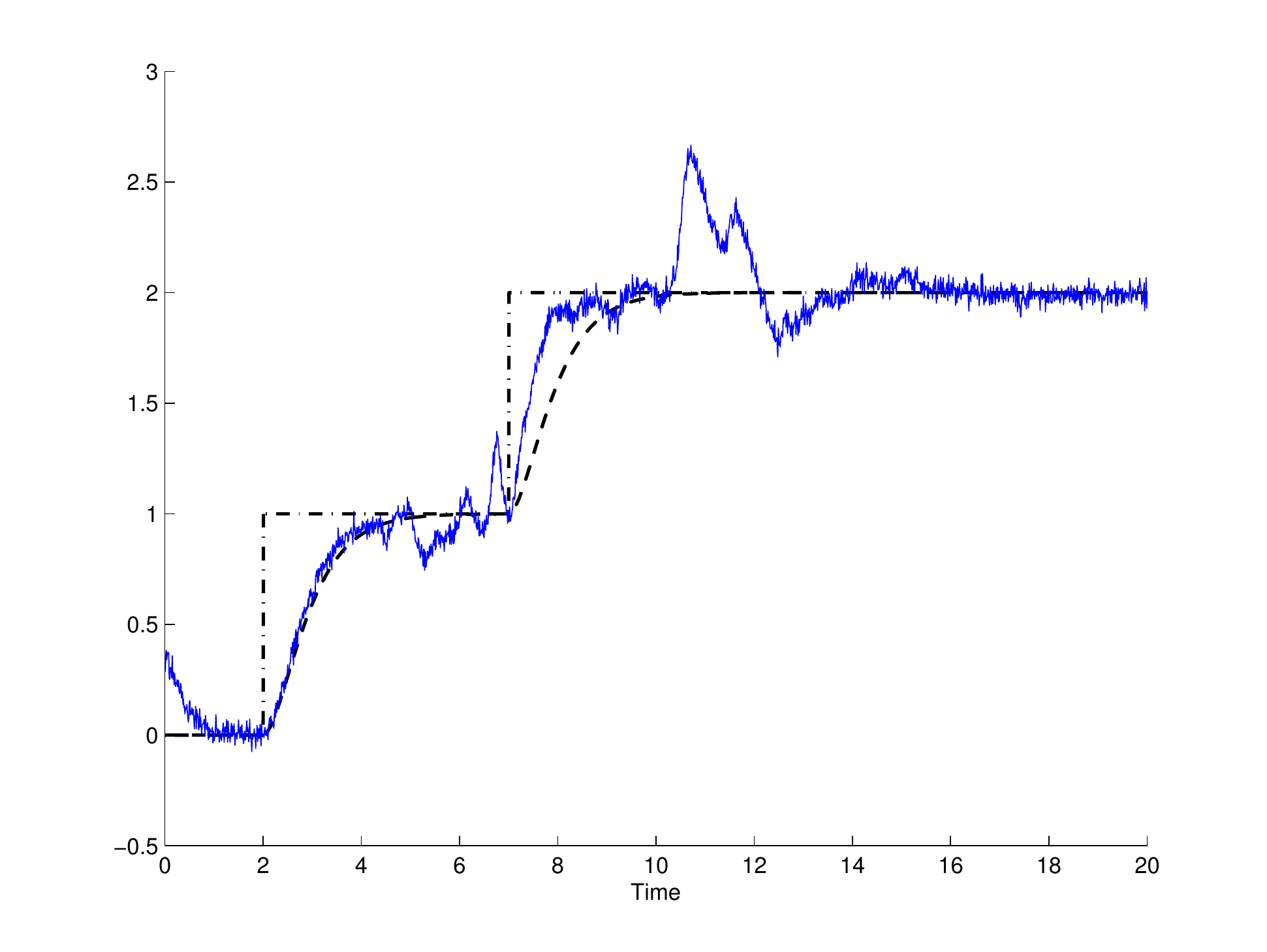}}}%
\caption{Delay system: Model-free control}%
\label{3X}
\end{center}
\end{figure}

Consider the system
$$\dot y(t)=y(t)+5y(t-\tau)+u$$
where moreover the delay $\tau$, $0\leq \tau \leq 5\text{s}$ is not
assumed to be
\begin{itemize}
\item known,
\item constant.
\end{itemize}
Set for the numerical simulations (see Figure \ref{3XD})
$$\tau (t)=\tau(t-T_e) +10T_e \text{sign}(N(t)), \quad \tau(0)=2.5\text{s}$$
where $N$ is a zero-mean Gaussian distribution with standard
deviation $1$.
An iP where $K_P = 1$ is deduced from the ultra-local model $\dot{y}
= F+u$. The results depicted in Figure \ref{3X} are quite
satisfactory.

\begin{remark}
Extending the above control strategy to non-linear systems and to
neutral systems is straightforward. It will
not be developed here. 
\end{remark}
\begin{remark}\label{practdelay}
The delay appearing with the hydro-electric power plants studied by \cite{edf} was taken into account via an empirical knowledge of the process.
Some numerical tabulations were employed in order to get in some sense ``rid'' of the delay. Such a viewpoint might be the most realistic one in industry.
\end{remark}
\begin{remark}
We only refer here to ``physical'' delays and
not to the familiar approximation in engineering of ``complex''
systems via delays ones (see, \textit{e.g.}, \cite{shin}). Let us emphasize that this type of 
approximation is loosing its importance in our setting.
\end{remark}

\subsection{A one-dimensional semi-linear heat equation}\label{fourier}
\begin{figure}
\begin{center}
\subfigure[Time evolution without measurement noise]{
\resizebox*{10cm}{!}{\includegraphics{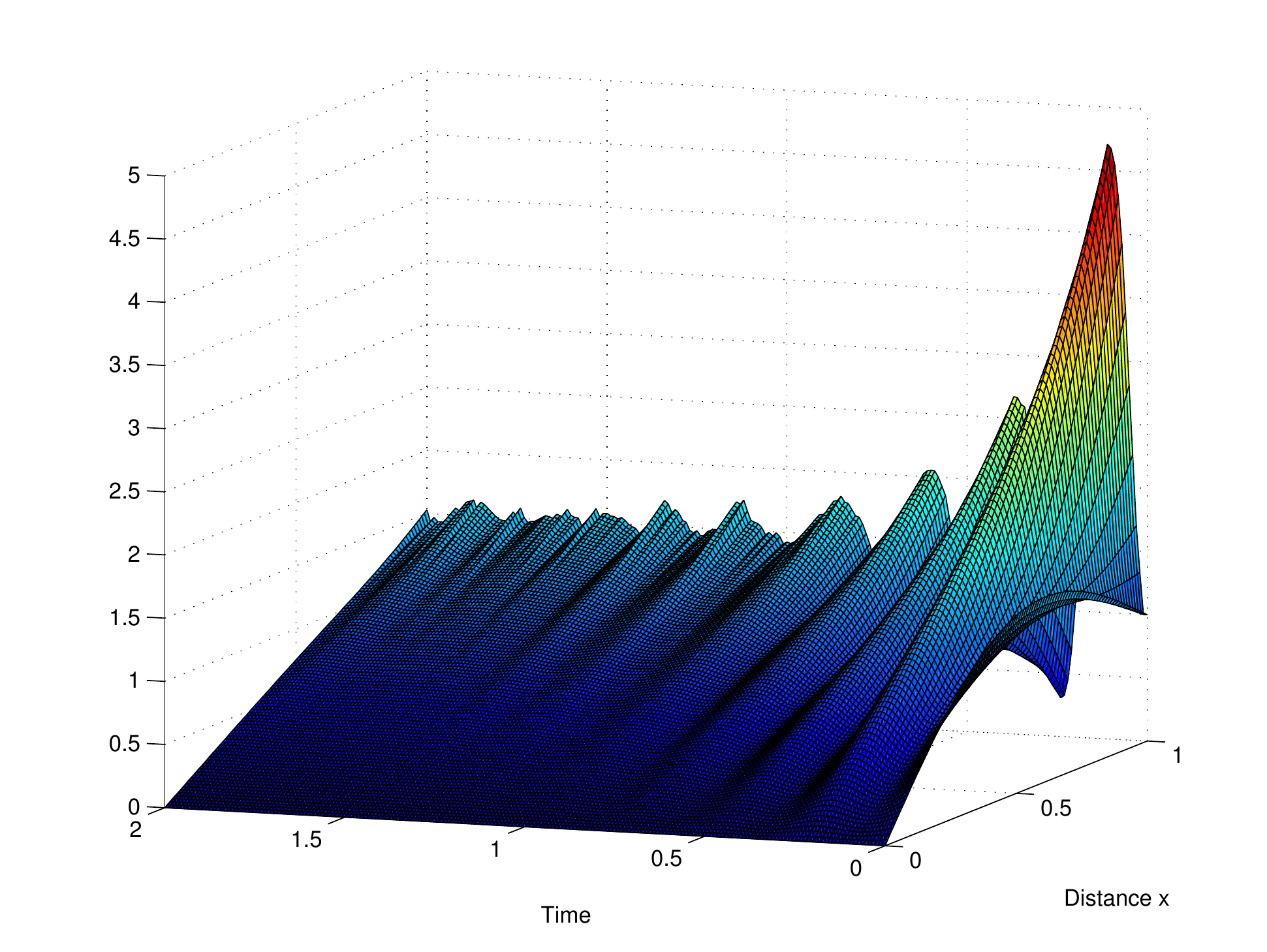}}}\\%
\subfigure[Control $u(t)$]{
\resizebox*{8cm}{!}{\includegraphics{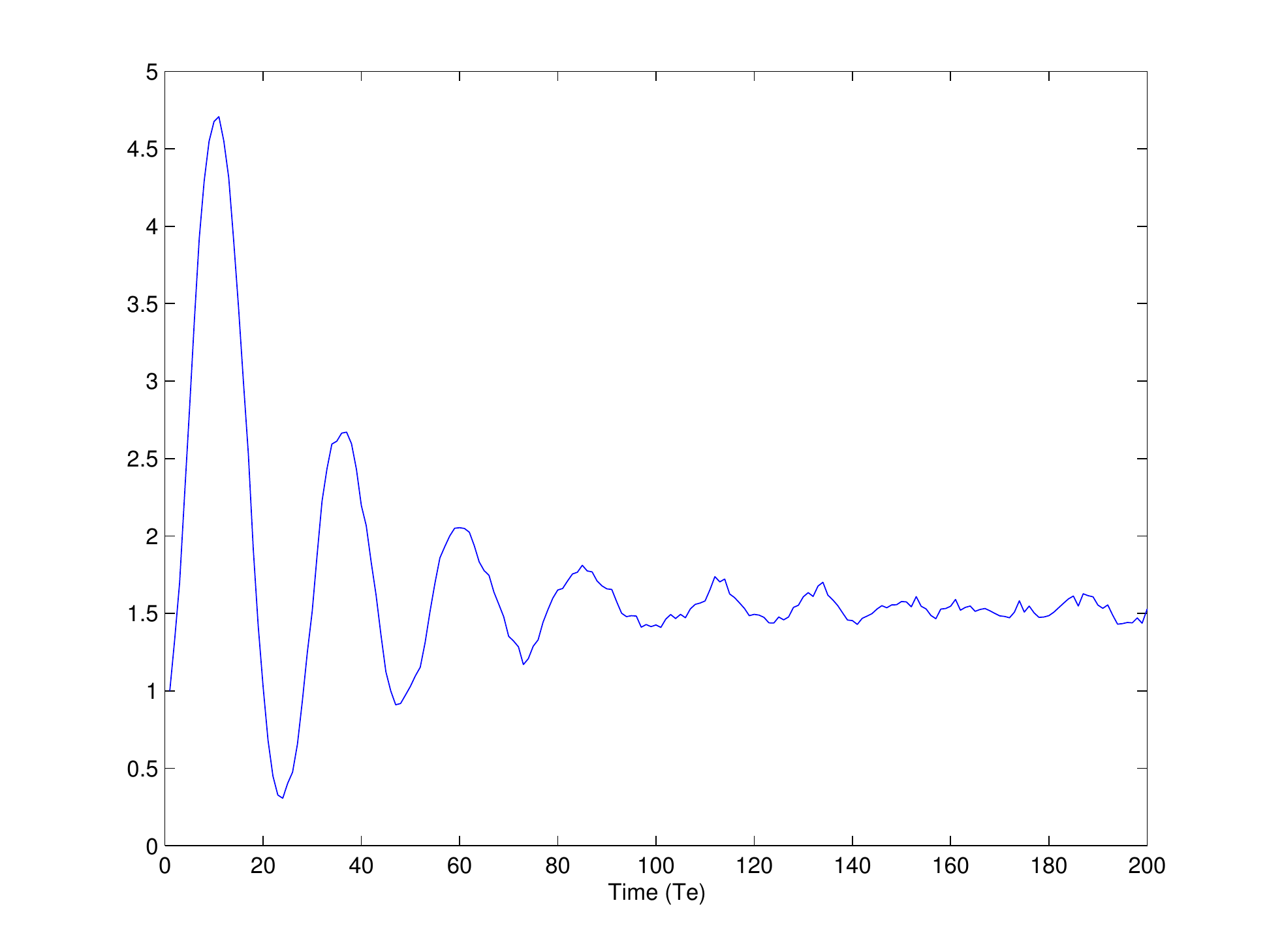}}}%
\subfigure[Controlled heat at distance $x_c$ (--, blue), setpoint (- ., black), and reference (- -, red)]{
\resizebox*{8cm}{!}{\includegraphics{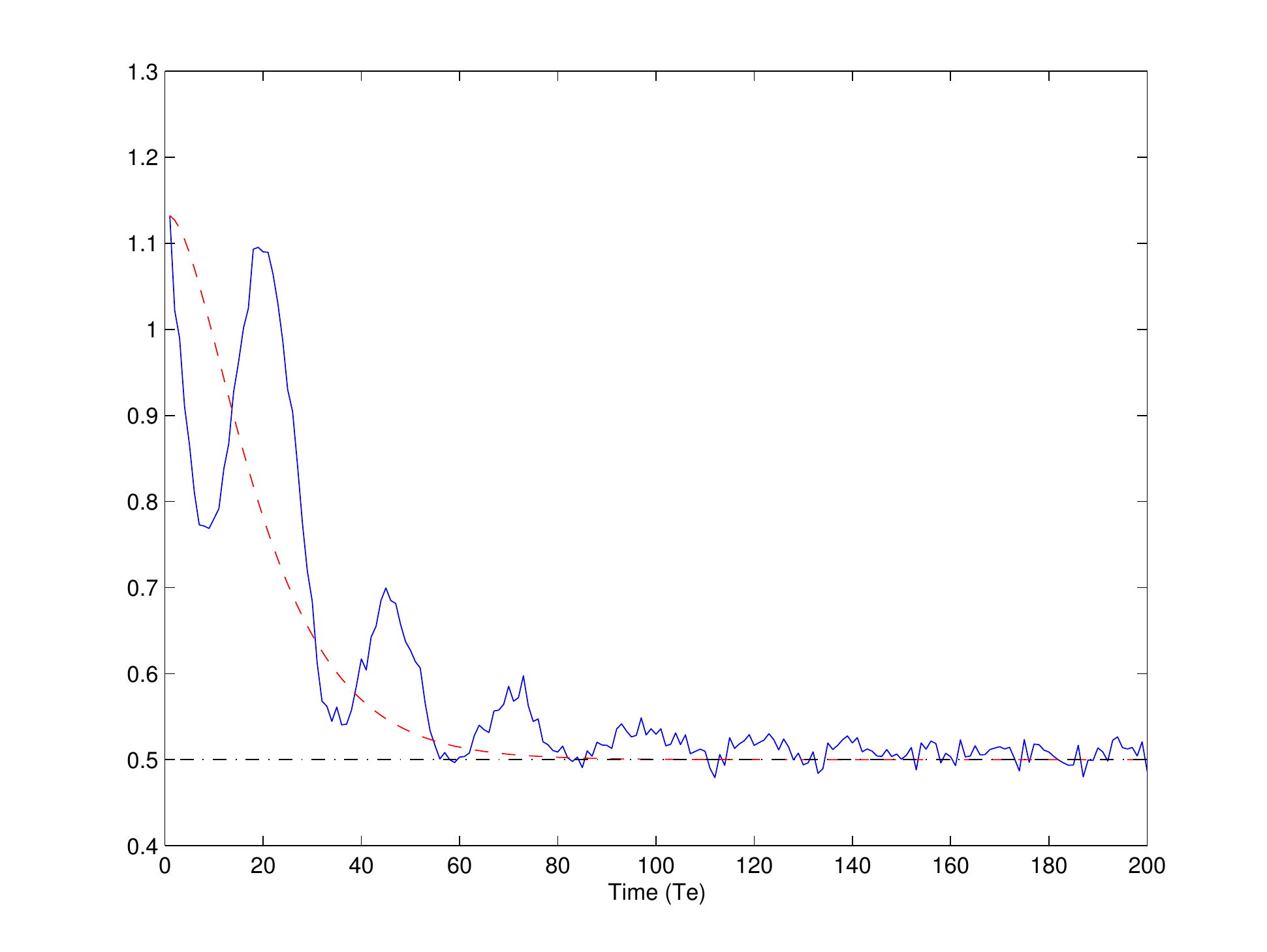}}}%
\caption{Heat equation: scenario 1}%
\label{pdes1}
\end{center}
\end{figure}

\begin{figure}
\begin{center}
\subfigure[Time evolution without measurement noise]{
\resizebox*{10cm}{!}{\includegraphics{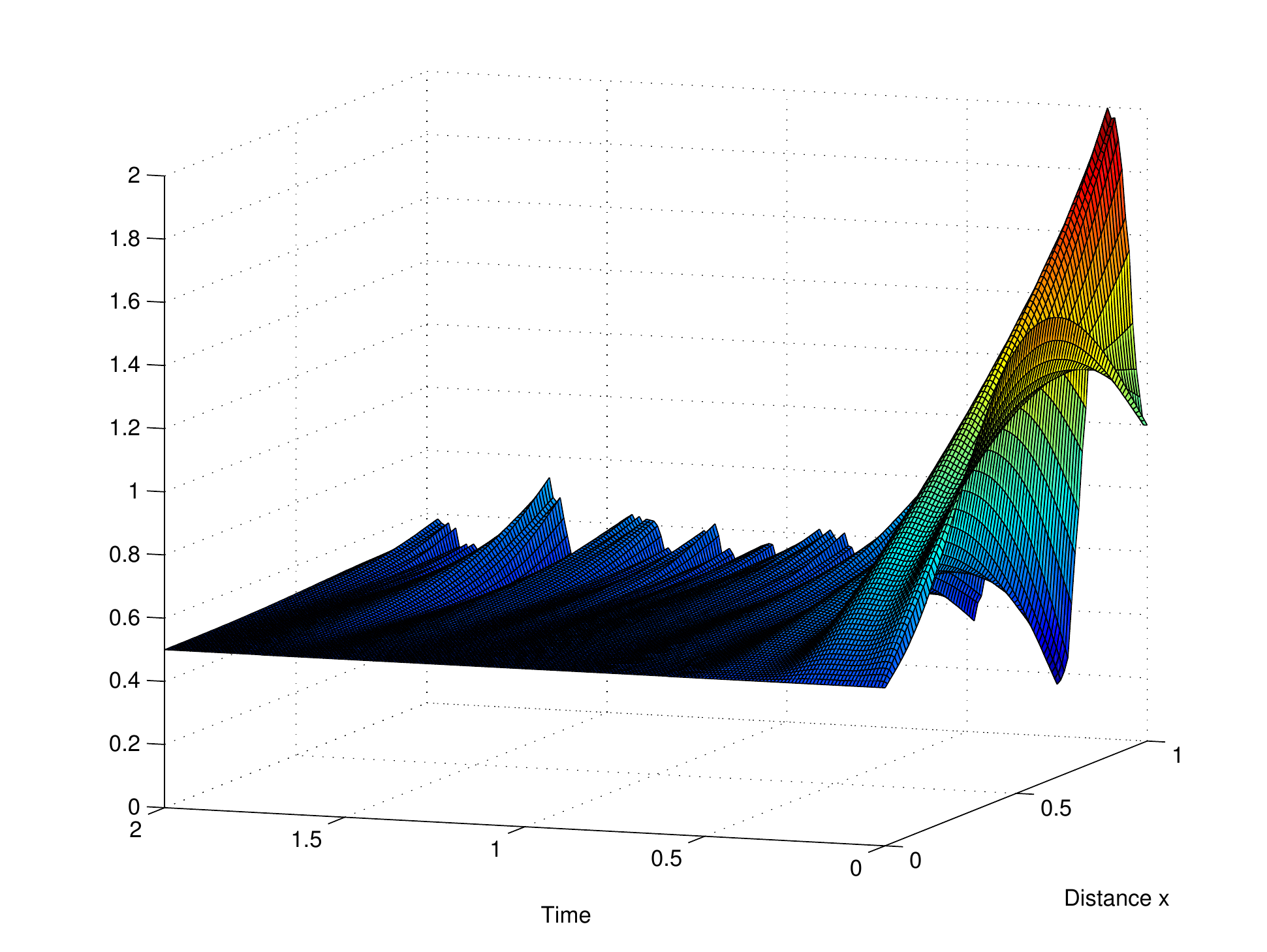}}}\\%
\subfigure[Control $u(t)$]{
\resizebox*{8cm}{!}{\includegraphics{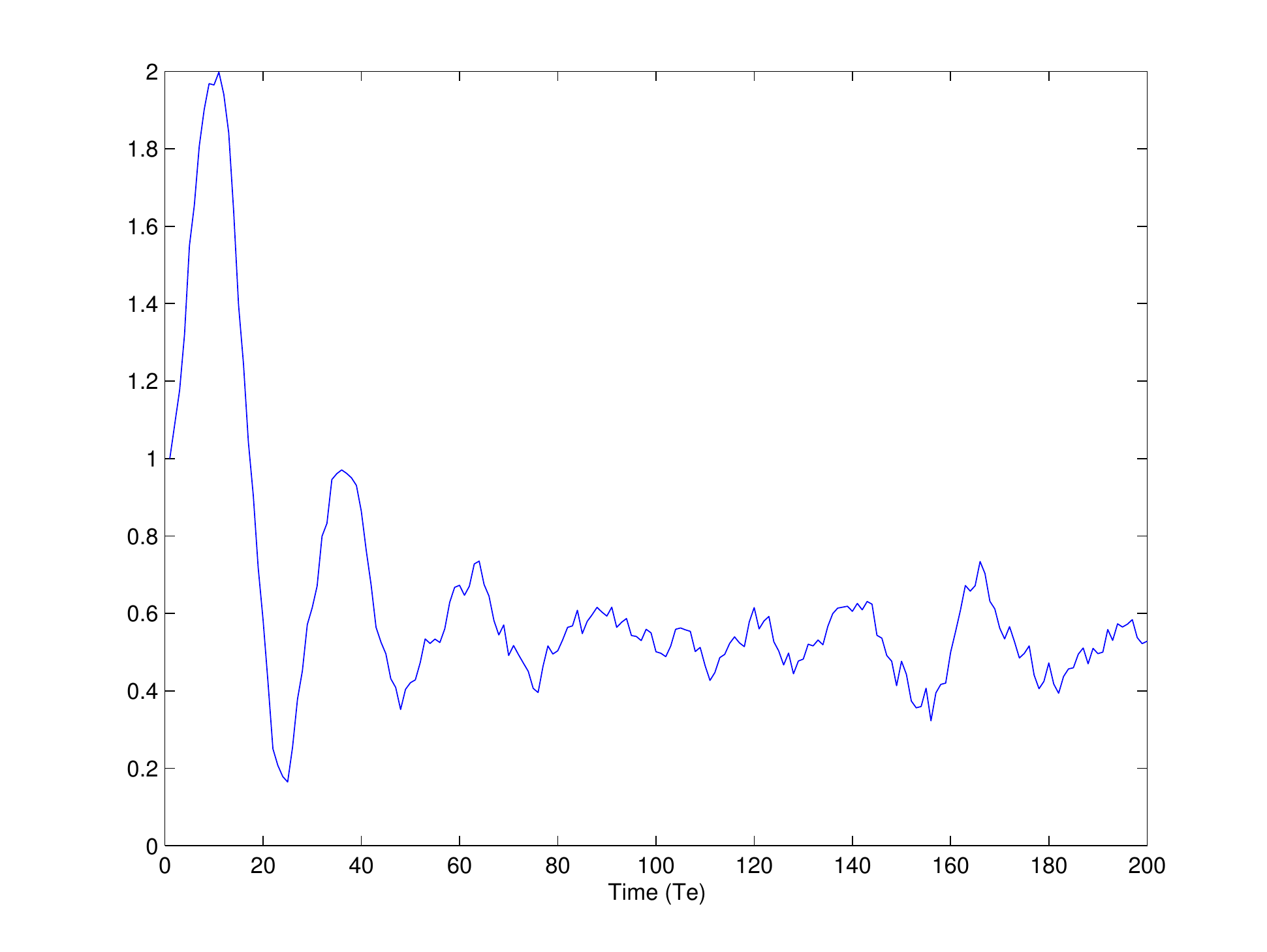}}}%
\subfigure[Controlled heat at distance $x_c$ (--, blue), setpoint (- ., black), and reference (- -, red)]{
\resizebox*{8cm}{!}{\includegraphics{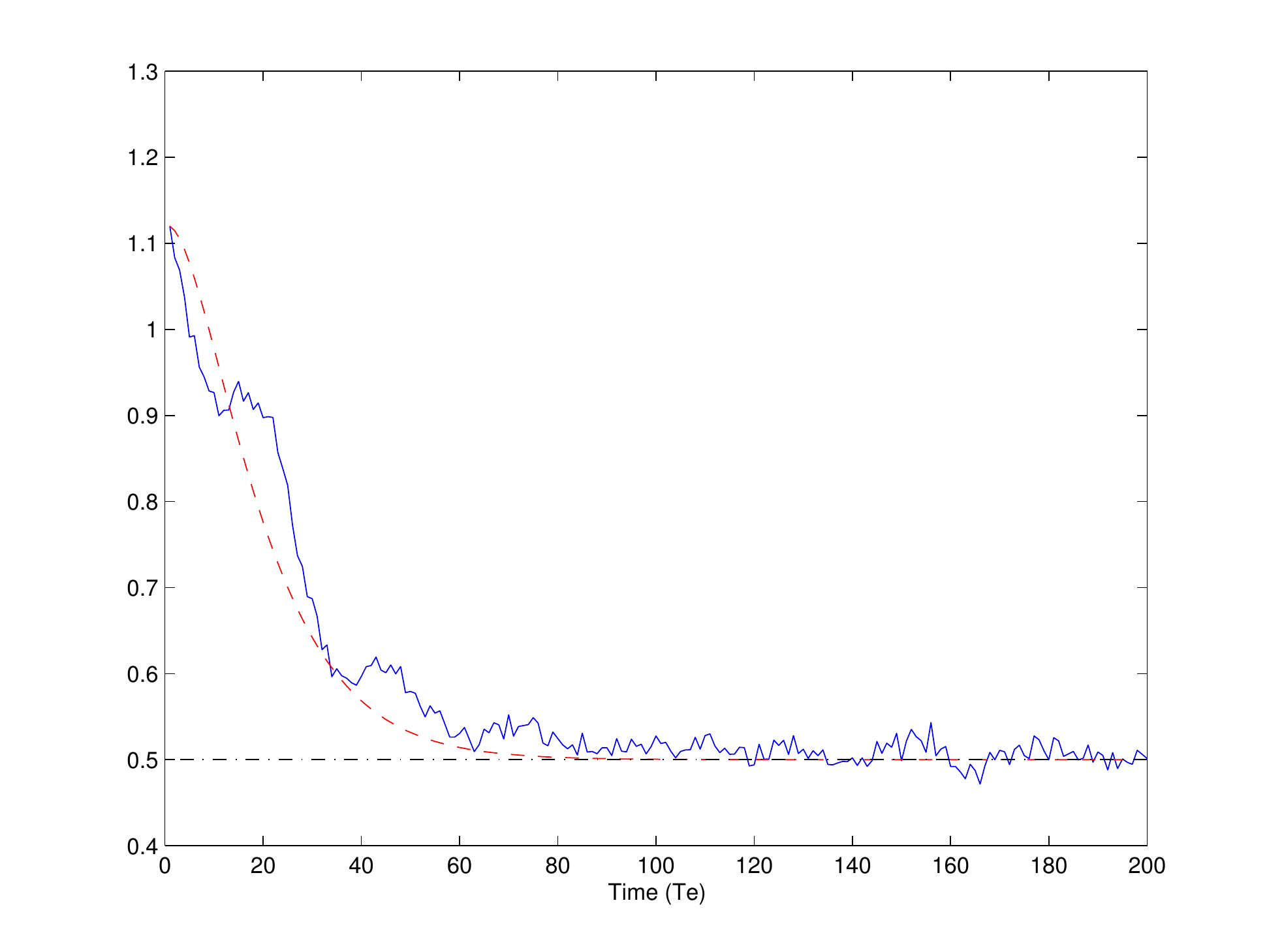}}}%
\caption{Heat equation: scenario 2}%
\label{pdes2}
\end{center}
\end{figure}

\begin{figure}
\begin{center}
\subfigure[Time evolution without measurement noise]{
\resizebox*{10cm}{!}{\includegraphics{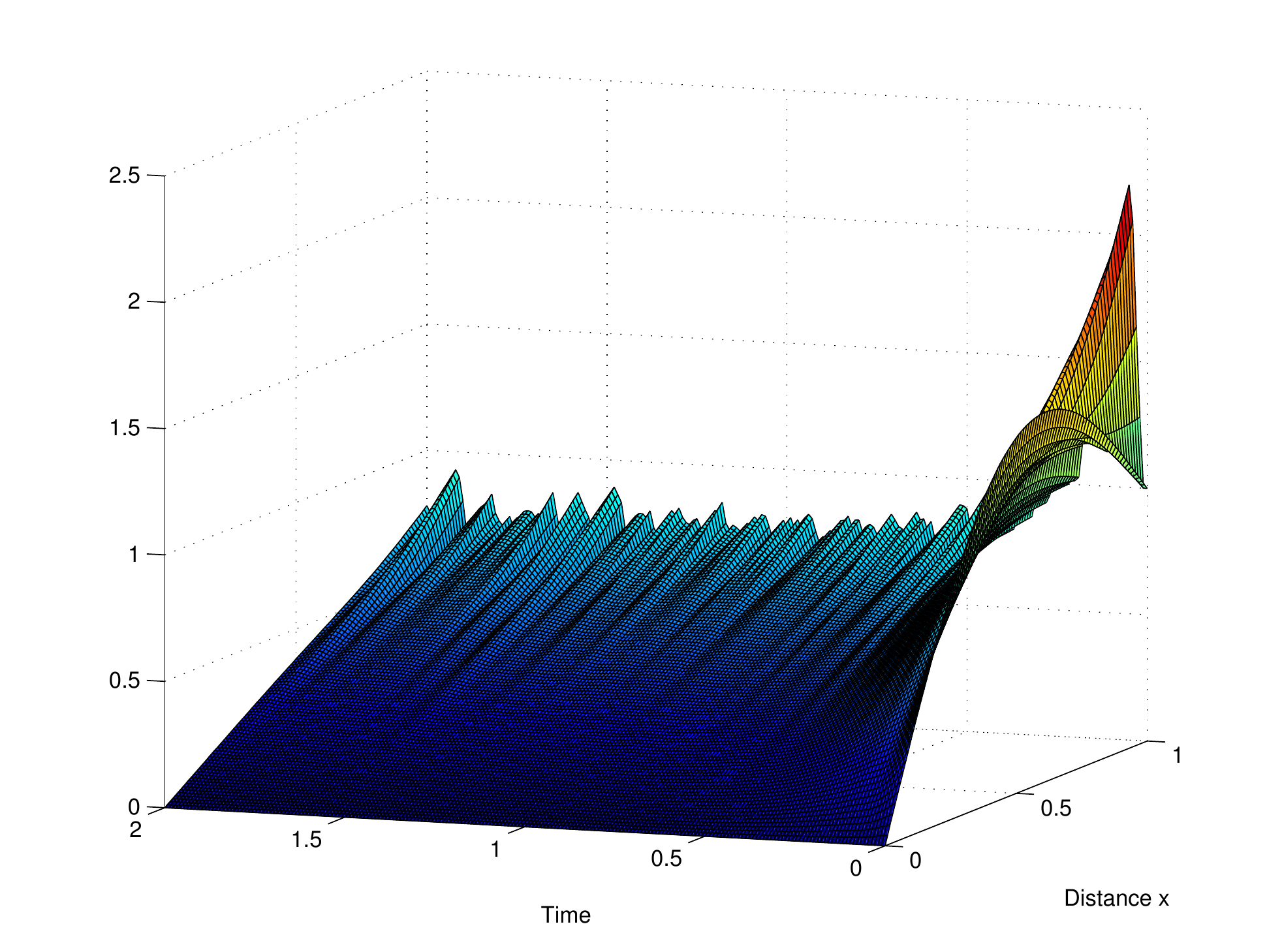}}}\\%
\subfigure[Control $u(t)$]{
\resizebox*{8cm}{!}{\includegraphics{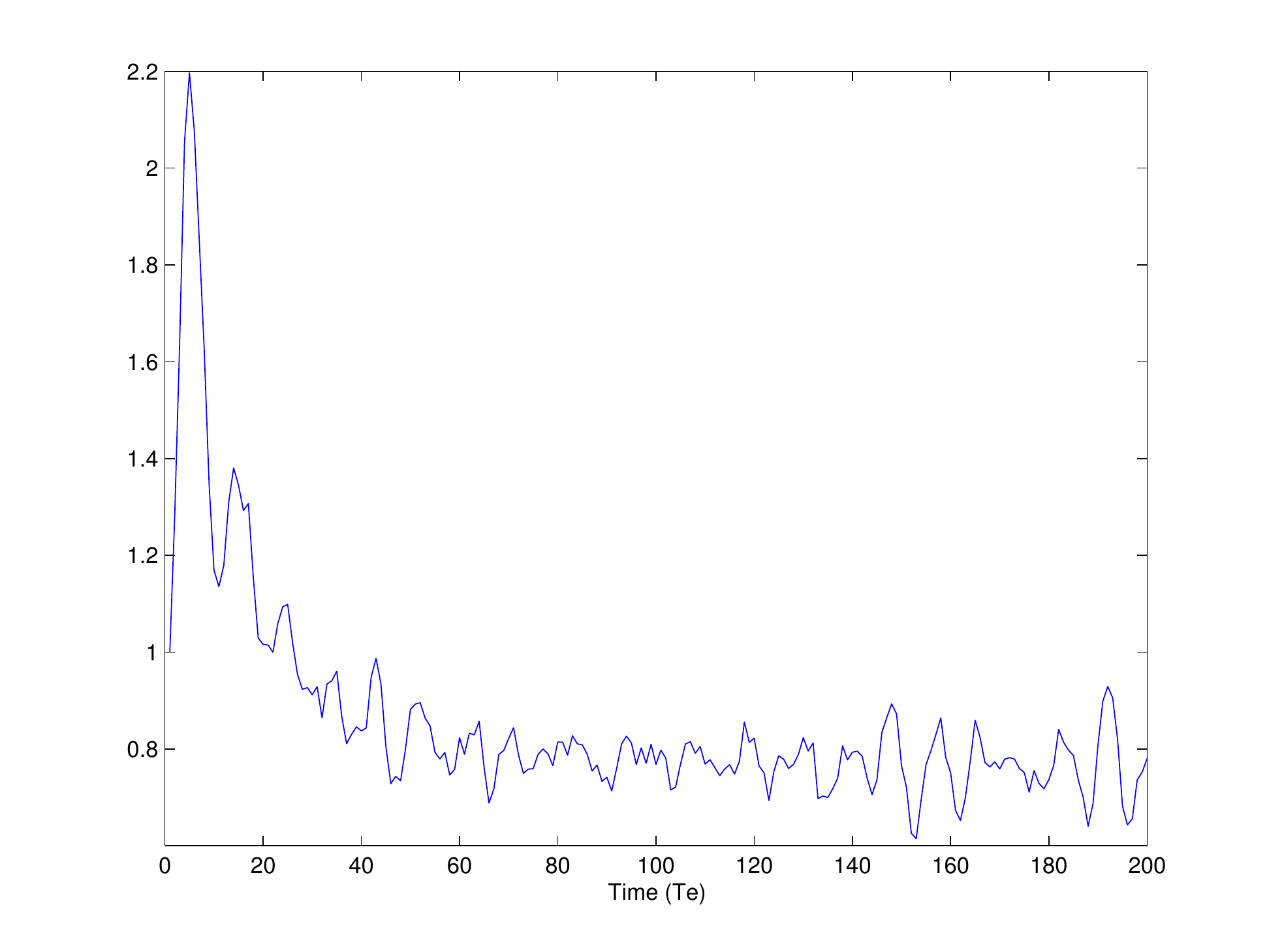}}}%
\subfigure[Controlled heat at distance $x_c$ (--, blue), setpoint (- .,black), and reference (- -, red)]{
\resizebox*{8cm}{!}{\includegraphics{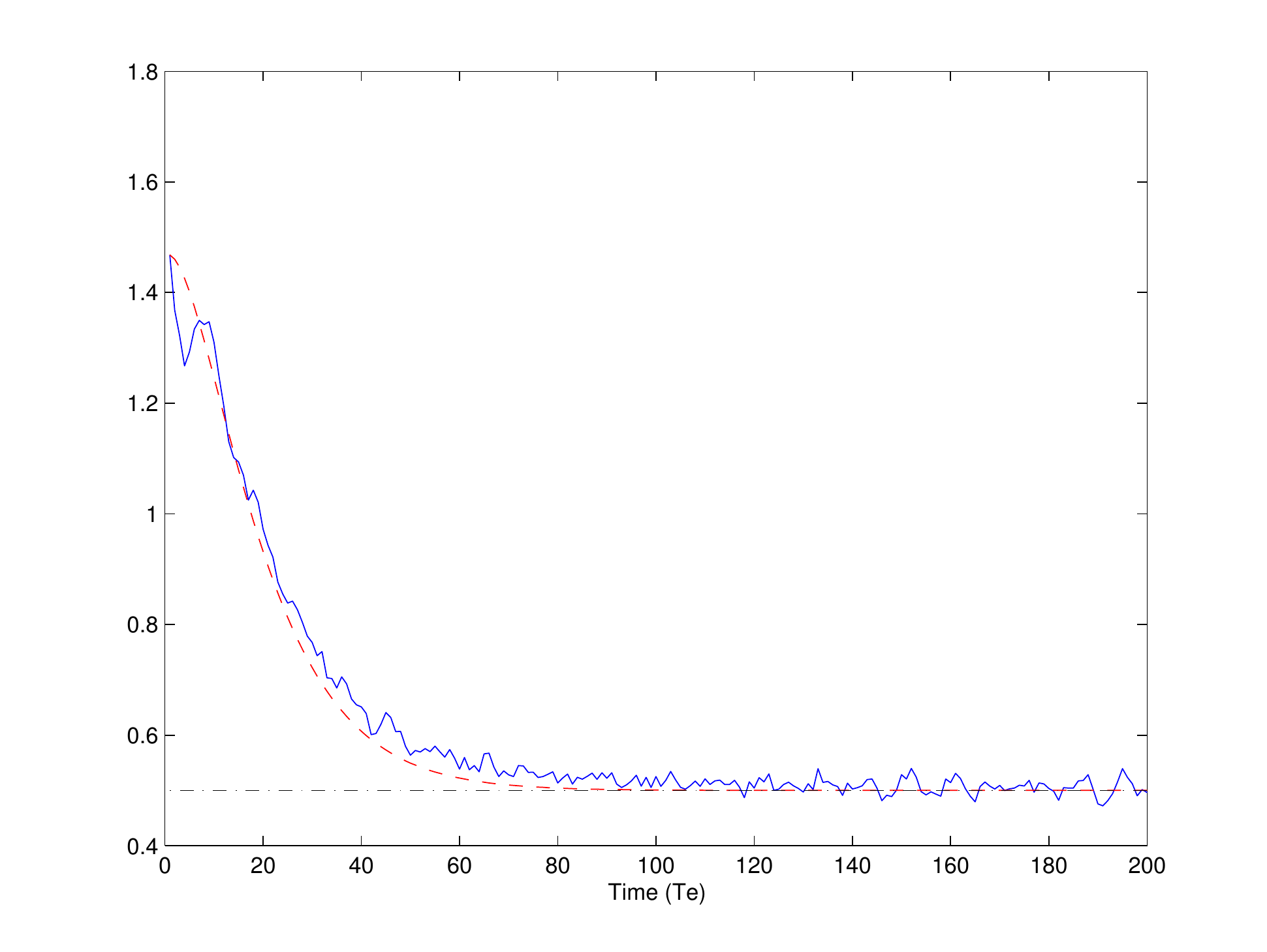}}}%
\caption{Heat equation: scenario 3}%
\label{pdes3}
\end{center}
\end{figure}

\begin{figure}
\begin{center}
\subfigure[Time evolution without measurement noise]{
\resizebox*{10cm}{!}{\includegraphics{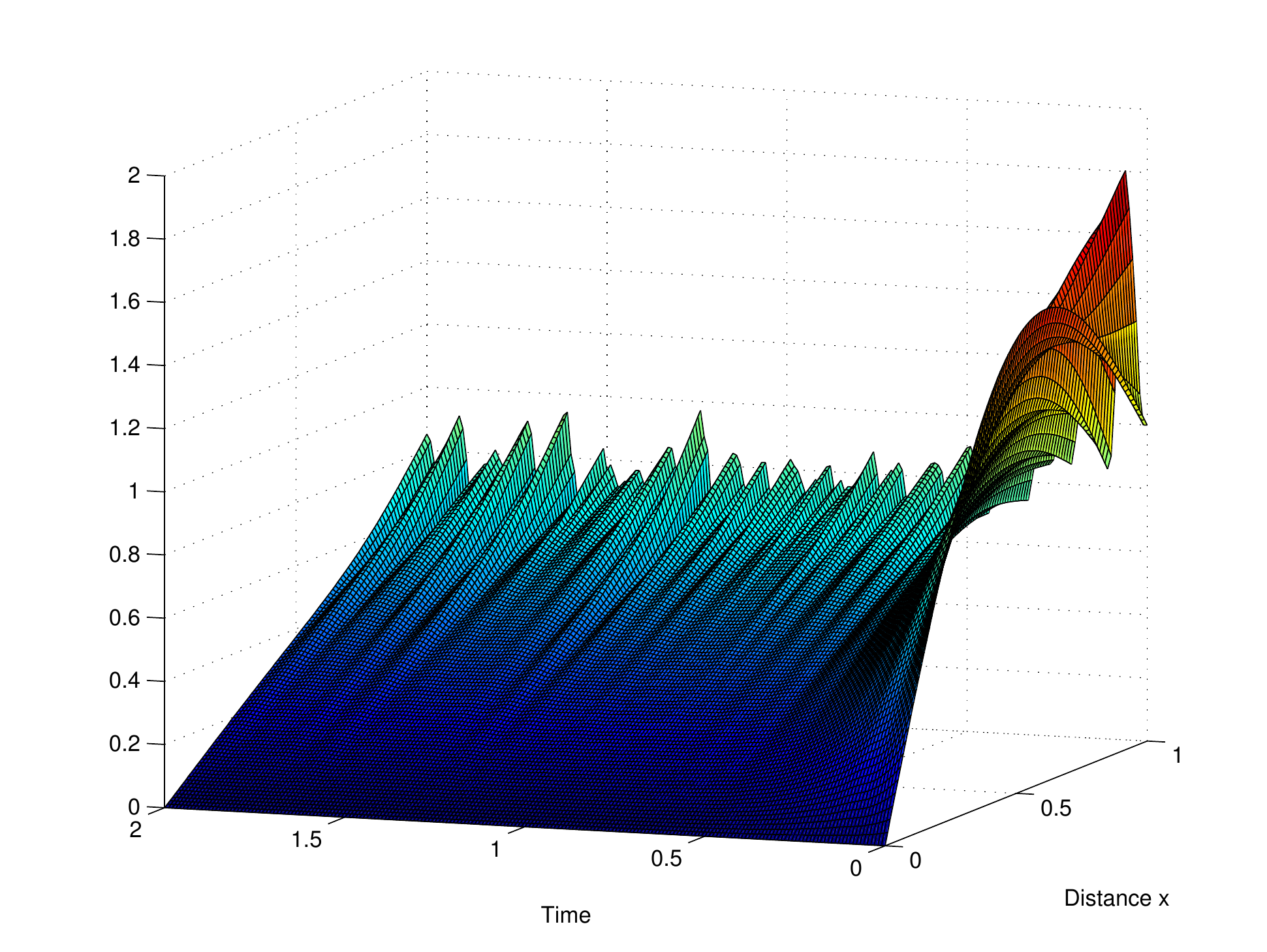}}}\\%
\subfigure[Control $u(t)$]{
\resizebox*{8cm}{!}{\includegraphics{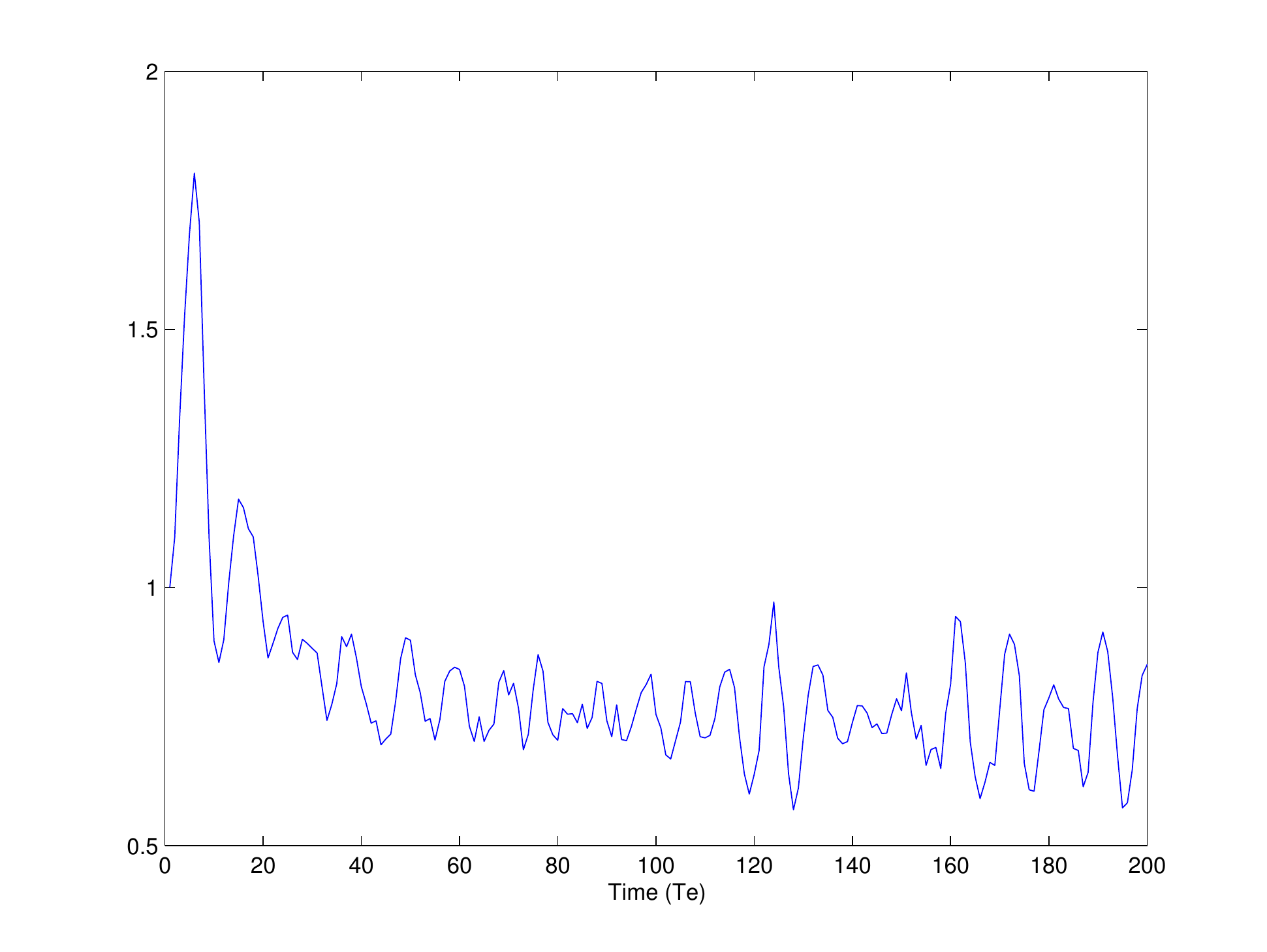}}}%
\subfigure[Controlled heat at distance $x_c$ (--, blue), setpoint (- ., black), and reference (- -, red)]{
\resizebox*{8cm}{!}{\includegraphics{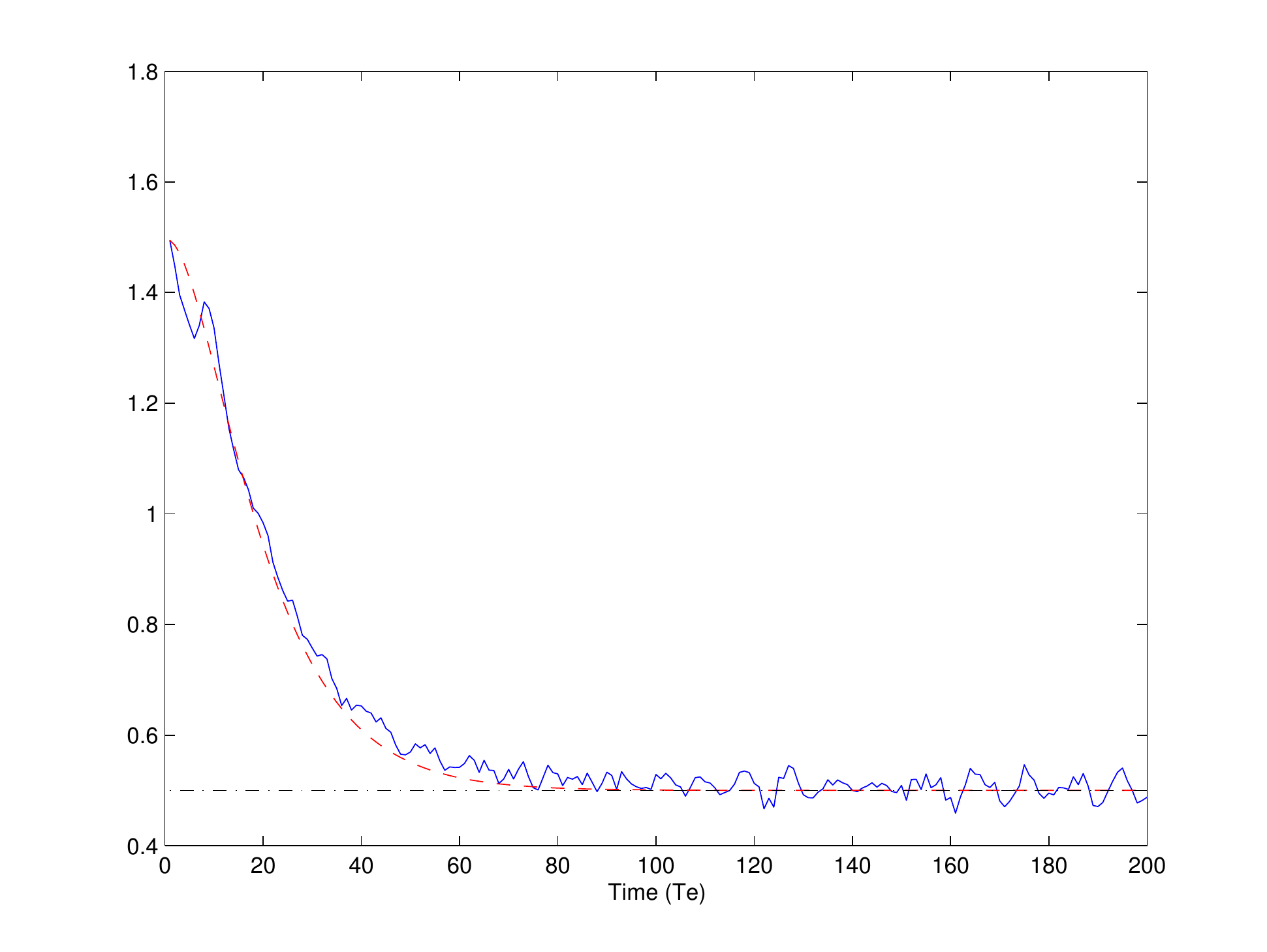}}}%
\caption{Heat equation: scenario 4}%
\label{pdes4}
\end{center}
\end{figure}

The heat equation is certainly one of the most studied topic in
mathematical physics. It would be pointless to review its
corresponding huge bibliography even in the control domain, where
many of the existing high-level control theories have been tested.
Consider with \cite{Coron} the one-dimensional semi-linear
heat equation
\begin{equation}
\frac{\partial w}{\partial t}=\frac{\partial^2 w}{\partial x^2}+f(w) 
\end{equation}
where
\begin{itemize}
\item $0 \leq x \leq L$,
\item $w(t,0)=c$, 
\item $w(t,L)=u(t)$ is the control variable,
\item $w(0,x)= sin(\pi x)+(u(0)-c)x+c$, where $c \in \mathbb{R}$ is a constant..
\end{itemize}
We want to obtain given
time-dependent temperature at $x = x_c$. Consider the following
scenarios:
\begin{enumerate}
\item $x_c=1/3 L$, $f=0$, $c=0$,
\item $x_c=1/3 L$, $f=0$, $c=0.5$,
\item $x_c=2/3 L$, $f=0$, $c=0$,
\item $x_c=2/3 L$, $f=y^3$, $c=0$.
\end{enumerate}
The control synthesis is achieved thanks to the elementary
one-dimensional ultra-local model
$$\dot y=F + 10 u$$
and the straightforward iP, where $K_p=10$. The four numerical simulations, displayed by Figures \ref{pdes1}, \ref{pdes2}, \ref{pdes3}, and \ref{pdes4}, are quite
convincing.

\subsection{A peculiar non-minimum phase system}\label{nm}
\begin{figure}
\begin{center}
\subfigure[Control]{
\resizebox*{15cm}{!}{\includegraphics{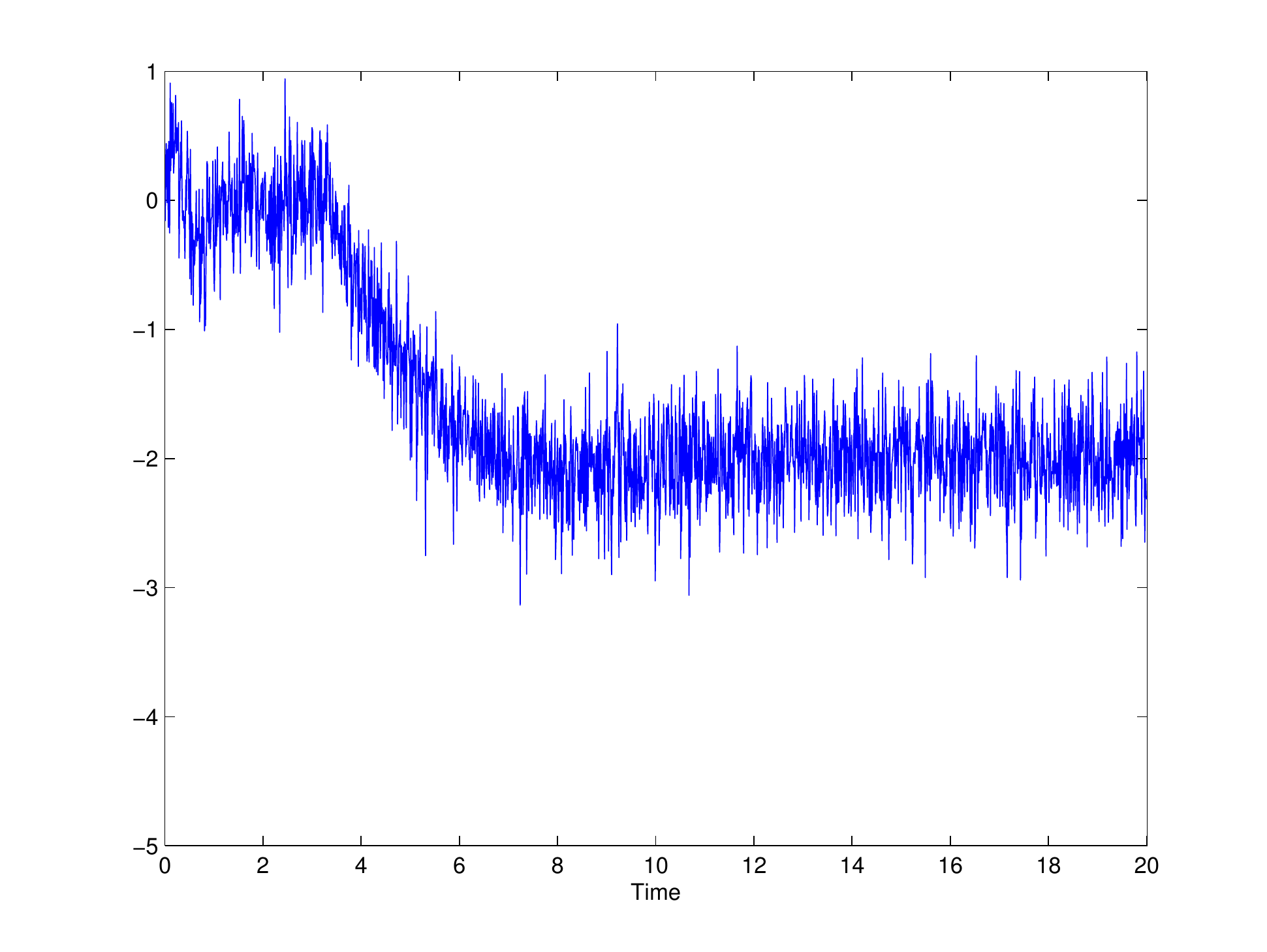}}}\\
\subfigure[Setpoint (- .,black) and Output (--,blue)]{
\resizebox*{15cm}{!}{\includegraphics{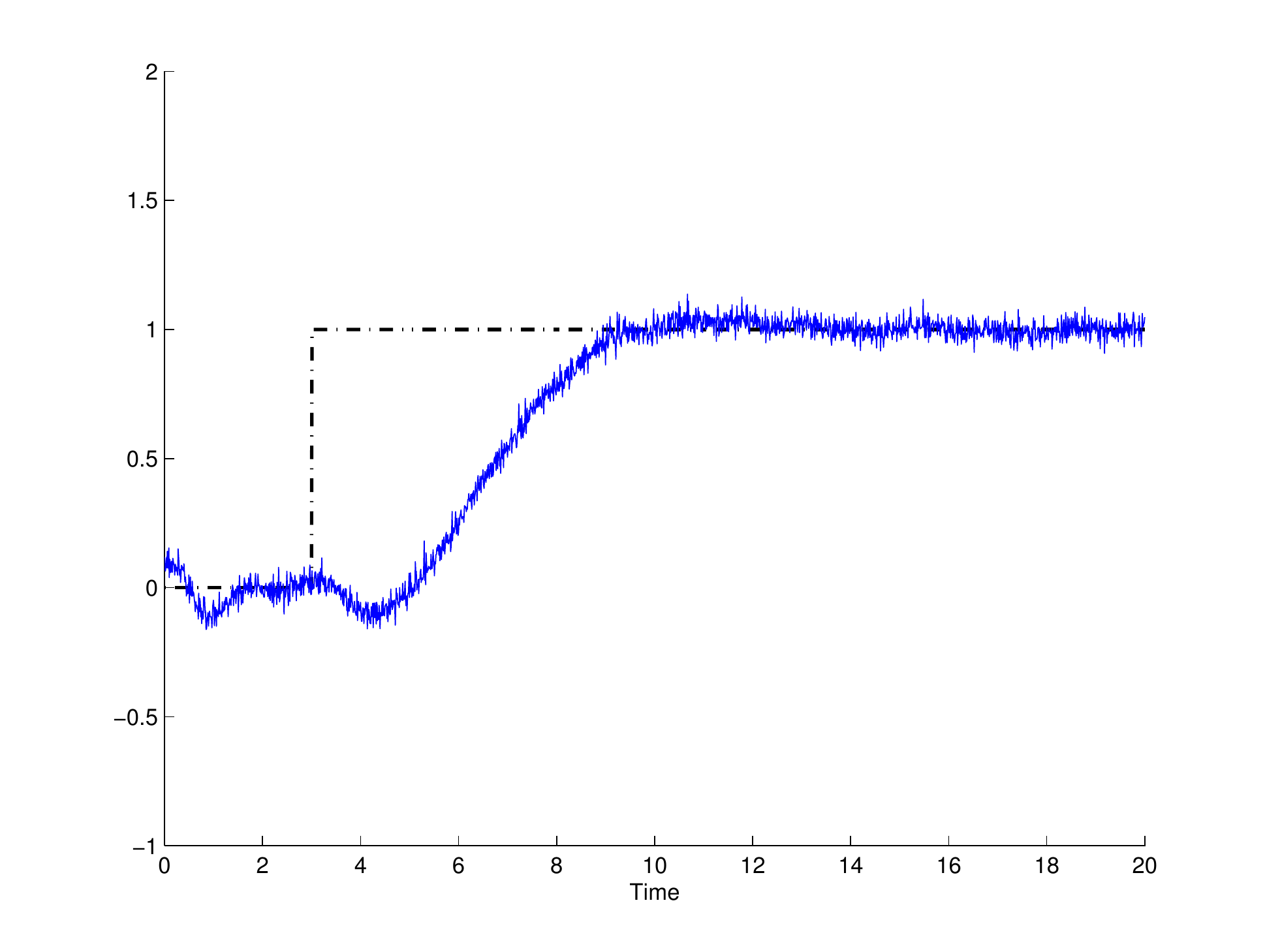}}}%
\caption{Non-minimum phase linear system: Model-free control}%
\label{7X}
\end{center}
\end{figure}

Consider the non-minimum phase system defined by the transfer function
\begin{equation}\label{transferPNM}
\frac{(s-1)}{(s+1)(s+2)}
\end{equation}
Utilize Equations \eqref{2i} and \eqref{i2}. Set $\alpha=-\beta=10$, $K_P=3$ and $K_I=K_{II}=5$. Figure \ref{7X} displays good performances.  

\begin{remark}
It is easy to check that the above calculations work only for a single unstable zero, like in Equation \eqref{transferPNM}. Our approach cannot be extended to arbitrary non-minimum phase systems.
\end{remark}

\begin{remark}\label{remnm}
It is well known that the control synthesis of a non-minimum phase system is even a difficult task with a perfectly known mathematical model. Among the many solutions which have been suggested in the literature, let us mention 
a flatness-based output change (see \cite{predict,sira3}). When a mathematical model is unknown or poorly known, the non-minimum phase character of an output cannot be 
deduced mathematically but only via a ``bad'' qualitative behavior of this output. Selecting a minimum phase output, \textit{i.e.}, an output with ``good'' qualitative properties, might be 
a more realistic alternative. It necessitates nevertheless an excellent ``practical'' knowledge of the plant behaviour.
\end{remark}

\section{Connections between classic and intelligent
controllers}\label{connect}
The results below connect classic PIDs to our intelligent controllers. They explain therefore why classic PIDs are used in rather arbitrary industrial situations thanks to a fine gain tuning, which 
might be quite difficult to achieve in practice. 
\subsection{PI and iP}\label{1}
\subsubsection{A crude sampling of PIs}
%
Consider the classic continuous-time PI controller
\begin{equation}\label{cpi}
  u (t) = k_p e(t) + k_i \int e(\tau) d\tau
\end{equation}
A crude sampling of the integral $\int e(\tau) d\tau$ through a
Riemann sum $I(t)$ leads to
$$
\int e(\tau) d\tau \simeq  I(t) = I(t-h) + h e(t)
$$
where $h$ is the sampling interval. The corresponding discrete form
of Equation \eqref{cpi} reads:
$$
u(t) = k_p e(t) + k_i I(t) = k_p e(t) + k_i I(t-h) + k_i h e(t)
$$
Combining the above equation with $$u(t-h) = k_p e(t-h) + k_i
I(t-h)$$ yields
\begin{equation}
\label{eqPIRiemannDiscrSix} u(t) = u(t - h) + k_p \left( e(t) - e(t
- h) \right) + k_i h e(t)
\end{equation}

\begin{remark}
A trivial sampling of the ``velocity form'' of Equation \eqref{cpi}
\begin{equation*}\label{cpid}
\dot{u} (t) = k_p \dot{e}(t) + k_i e(t)
\end{equation*}
yields
$$
\dfrac{u(t) - u(t - h)}{h} =  k_p  \left(\dfrac{e(t) - e(t -
h)}{h}\right) + k_i  e(t)
$$
which is equivalent to Equation \eqref{eqPIRiemannDiscrSix}.
\end{remark}

\subsubsection{Sampling iPs} Utilize, if $\nu = 1$, the iP, which may be rewritten as
$$
u (t) = \frac{{\dot y}^\ast (t) - F + K_P e(t)}{\alpha}
$$
Replace $F$ by ${\dot y}(t) - \alpha u (t-h)$ and therefore by
$$\frac{y(t) - y(t-h)}{h} - \alpha u (t-h)$$
It yields
\begin{equation}
\label{eqDiscr_i-POne} u (t) = u (t - h) - \frac{e(t) -
e(t-h)}{h\alpha} + \dfrac{K_P}{\alpha}\, e(t)
\end{equation}
\subsubsection{Comparison}
{\bf FACT}.- Equations \eqref{eqPIRiemannDiscrSix} and
\eqref{eqDiscr_i-POne} become {\bf identical} if we set
\begin{align}
\label{eqPI_i-P_corresp} k_p &= - \dfrac{1}{\alpha h}, \quad k_i =
\dfrac{K_P}{\alpha h}
\end{align}

\begin{remark}
It should be emphasized that the above property, defined by
Equations \eqref{eqPI_i-P_corresp}, does not hold for
continuous-time PIs and iPs. This equivalence is strictly related to
time sampling, {\em i.e.}, to computer implementation, as
demonstrated by taking $h \downarrow 0$ in Equations
\eqref{eqPI_i-P_corresp}.
\end{remark}

\subsection{PID and iPD}
Extending the calculations of Section \ref{1} is quite obvious. The
velocity form of the PID
$$u(t) = k_p e(t) + k_i \int e(\tau) d\tau + k_d \dot{e}$$ reads
$\dot u(t) = k_p {\dot e}(t) + k_i e(t) + k_d {\ddot e}(t)$. It
yields the obvious sampling
\begin{equation}
\label{eqDiscr_PID}
  u(t) = u(t - h) + k_p h {\dot e} (t) + k_i h e(t) +
            k_d h {\ddot e}(t)
\end{equation}
$\nu = 2$ on the other hand, Equation yields $u (t) =
\dfrac{1}{\alpha} \left( {\ddot y}^\ast (t) - F + K_P e(t) + K_D
{\dot e}(t) \right)$. From the computer implementation $F =
\ddot{y}(t) - \alpha u(t - h)$, we derive
\begin{equation}\label{iPDd}
u (t) = u (t - h) - \dfrac{1}{\alpha} {\ddot e}(t)  +
\dfrac{K_P}{\alpha} e(t) + \dfrac{K_D}{\alpha} {\dot e}(t)
\end{equation}

\noindent{\bf FACT}.- Equations \eqref{eqDiscr_PID} and \eqref{iPDd}
become {\bf identical} if we set
\begin{equation}
\label{eqPID_i-PD_corresp} k_p = \dfrac{K_D}{\alpha h}, \quad k_i =
\dfrac{K_P}{\alpha h}, \quad k_d = - \dfrac{1}{\alpha h}
\end{equation}

\subsection{iPI and iPID}\label{2int}
\label{Seci-PIsi-PIDs} Equation \eqref{iPDd} becomes with the iPID
\begin{equation}\label{iPIDd}
u (t) = u (t - h) - \dfrac{1}{\alpha} {\ddot e}(t)  +
\dfrac{K_P}{\alpha} e(t) + \dfrac{K_I}{\alpha} \int e +
\dfrac{K_D}{\alpha} {\dot e}(t)
\end{equation}
Introduce the PI$^2$D controller
$$u (t) = k_p e (t) + k_i \int e(\tau)d\tau  + k_{ii}
\int\!\!\!\!\int e d\tau d\sigma + k_d \dot{e}(t)$$ The double
integral, which appears there, seems to be quite
uncommon in control engineering. To its velocity form ${\dot u} (t) =
k_p \dot{e} (t) + k_i e + k_{ii} \int e d\tau + k_d \ddot{e}(t)$
corresponds the sampling
$$
u (t) = u(t-h) + k_p h \dot{e} (t) + k_i h e + k_{ii} h \int e d\tau
+ k_d h \ddot{e}(t) $$ which is identical to Equation \eqref{iPIDd}
if one sets
\begin{align}
    \label{eqPID_i-PD_correspTwo}
    k_p &= \dfrac{K_D}{\alpha h}, \quad k_i = \dfrac{K_P}{\alpha h}, \quad
    k_{ii} = \dfrac{K_I}{\alpha h}, \quad k_d = - \dfrac{1}{\alpha h}
  \end{align}
The connection between iPIs and PI$^2$s follows at once.
\subsection{Table of correspondence} \label{SecCorrespondance} The
previous calculations yield the following correspondence (Table \ref{tblGainsCorresp})
between the gains of our various controllers: \ifdeuxcols
\begin{footnotesize}
\fi
\begin{table*}[htb]\label{table}
\begin{center}
\begin{tabular}{p{5ex}p{3ex}p{7ex}p{7ex}p{7ex}p{7ex}}
\toprule
         &       & iP            & iPD           & iPI & iPID          \\
\midrule
PI       & $k_p$ & $- 1/\alpha h$   &                &      &                \\
         & $k_i$ & $K_P/\alpha h$ &                &      &                \\
\midrule
PID      & $k_p$ &                & $K_D/\alpha h$   &      &                \\
         & $k_i$ &                & $K_P/\alpha h$ &      &                \\
         & $k_d$ &                & $- 1/\alpha h$ &      &                \\
\midrule
PI$^2$  & $k_p$    &             &                & $- 1/\alpha h$   &      \\
         & $k_i$    &             &                & $K_P/\alpha h$ &      \\
         & $k_{ii}$ &             &                & $K_I/\alpha h$ &      \\
\midrule
PI$^2$D & $k_p$    &             &                &      & $K_D/\alpha h$   \\
         & $k_i$    &             &                &      & $K_P/\alpha h$ \\
         & $k_{ii}$ &             &                &      & $K_I/\alpha h$ \\
         & $k_d$    &             &                &      & $- 1/\alpha h$ \\
\bottomrule
\end{tabular}
\end{center}
\caption{\label{tblGainsCorresp}Correspondence between the gains of
sampled classic and intelligent controllers.}
\end{table*}
\ifdeuxcols
\end{footnotesize}
\fi

\begin{remark}
Due to the form of Equation \eqref{cpi}, it should be noticed that the
tuning gains of the classic regulators ought to be negative.
\end{remark}

\section{Conclusion}\label{conclusion}

Several theoretical questions remain of course open. Let us mention some of them, which appear today to be most important:
\begin{itemize}
\item  The fact that multivariable systems were not studied here is due to a lack until now of concrete case-studies.
They should therefore be examined more closely.
\item Even if some delay and/or non-minimum phase examples were already successfully treated (see Sections \ref{del}, \ref{nm},  and (\cite{acc,edf,cifa-riachy})), a general understanding is still missing, like, 
to the best of our knowledge, with any other recent setting (see, \textit{e.g.}, \cite{astrom1,od}, and \cite{xu}). We believe as advocated in Remarks \ref{practdelay} and \ref{remnm} that 
\begin{itemize}
\item looking for a purely mathematical solution might be misleading,
\item taking advantage on the other hand of a ``good'' empirical understanding of the plant might lead to a more realistic track. 
\end{itemize}
\end{itemize}
It goes without saying that comparisons with existing approaches should be further explored. It has already been done with
\begin{itemize}
\item classic PIDs here, and by  \cite{brest,mil}  for some active spring and vehicles,
\item some aspects of \emph{sliding modes} by \cite{milan-sliding},
\item \emph{fuzzy control} for some vehicles by \cite{mil,vil2}.
\end{itemize}
Those comparisons were until now always favourable to our setting.

If model-free control and the corresponding intelligent controllers are further reinforced, especially by 
numerous fruitful applications, the consequences on the future development and teaching (see, \textit{e.g.}, the excellent textbook by \cite{murray}) of control theory might be dramatic:
\begin{itemize}
\item Questions on the structure and on the parameter identification of linear and nonlinear systems might loose their importance if the need of a ``good'' mathematical modeling is diminishing. 
\item Many effort on robustness issues with respect to a ``poor'' modeling and/or to disturbances may be viewed as obsolete and therefore less important. As a matter of fact those issues disappear to a large extent 
thanks to the continuously updated numerical values of  $F$ in Equation \eqref{ultralocal}. 
\end{itemize}
Another question, which was already raised by \cite{sofia}, should be emphasized.  Our model-free control strategy yields a 
straightforward regulation of industrial plants whereas the corresponding digital simulations need a reasonably accurate mathematical model in order to feed the computers. Advanced 
parameter identification and numerical analysis techniques might then be necessary tools  (see, \textit{e.g.}, \cite{edf,sofia}). This dichotomy between elementary control implementations and intricate computer
simulations  seems to the best of our knowledge to have been ignored until today. It should certainly be further dissected as a fundamental epistemological matter in engineering and, perhaps also, in other 
fields.

\label{lastpage}

\appendix
\section{An approximation property}\label{appendix}
\subsection{Functionals}\label{app1}
We restrict ourselves to a SISO system, \textit{i.e.}, to a system
with a single control variable $u$ and a single output $y$. Even
without knowing any ``good'' mathematical model we may assume that
the system corresponds to a \emph{causal}, or
\emph{non-anticipative}, \emph{functional}, \textit{i.e.}, for any
time instant $t
> 0$,
\begin{equation}\label{functional}
y(t) = \mathcal{F}\left( u(\tau) ~ | ~ 0 \leq \tau \leq t \right)
\end{equation}
where $\mathcal{F}$ depends on
\begin{itemize}
\item the past and the present, and not on the future,
\item various perturbations,
\item initial conditions at $t = 0$.
\end{itemize}

\begin{example}
A popular representation of rather arbitrary nonlinear systems in
engineering is provided by Volterra series (see, \textit{e.g.},
\cite{barrett}, \cite{rugh} and \cite{fl2}). Such a series may be
defined by
\begin{align*}\label{vs}
y(t) = & h_0(t) + \int_0^t h_1(t, \tau) u(\tau) d\tau + \\ &
\int_0^t \int_0^t h_2(t, \tau_2, \tau_1) u(\tau_2) u(\tau_1) d\tau_2
d\tau_1 + \dots \\
& \int_0^t \dots \int_0^t h_\nu(t, \tau_\nu, \dots \tau_1)
u(\tau_\nu) \dots u(\tau_1) d\tau_\nu \dots d\tau_1  \\ & + \dots
\end{align*}
Solutions of quite arbitrary differential equations may be expressed
as Volterra series.
\end{example}

\subsection{The Stone-Weierstra{\ss} theorem} Let
\begin{itemize}
\item $\mathcal{I} \subset [ 0, + \infty [$ be a compact subset,
\item $\mathcal{C} \subset C^0 (\mathcal{I})$ be a compact subset,
where $C^0 (\mathcal{I})$ is the space of continuous functions
$\mathcal{I} \rightarrow \mathbb{R}$, which is equipped with the
topology of uniform convergence.
\end{itemize}
Consider the Banach $\mathbb{R}$-algebra $\mathfrak{S}$ of
continuous causal functionals \eqref{functional} $\mathcal{I} \times
\mathcal{C} \rightarrow \mathbb{R}$. If a subalgebra contains a
non-zero constant element and separates points in $\mathcal{I}
\times \mathcal{C}$, then it is dense in $\mathfrak{S}$ according to the Stone-Weierstra{\ss} theorem (see,
\textit{e.g.}, the excellent textbooks by \cite{choquet} and \cite{rudin2}).

\subsection{Algebraic differential equations} Let $\mathfrak{A}
\subset \mathfrak{S}$ be the set of functionals which satisfy an
algebraic differential equation of the form
\begin{equation}\label{eq}
E(y, \dot{y}, \dots, y^{(a)}, u, \dot{u}, \dots, u^{(b)}) = 0
\end{equation}
where $E$ is a polynomial function of its arguments with real
coefficients. Satisfying Equation \eqref{eq} is equivalent saying
that $y$ is \emph{differential algebraic} over the
\emph{differential field} $\mathbb{R} \langle u \rangle$. 
\begin{remark}
Remind that a
\emph{differential field} (see, \textit{e.g.}, the two following books by \cite{cl} and \cite{kolchin}, 
and the papers by \cite{delaleau}, \cite{nl} and \cite{flmr}) is a commutative field which is
equipped with a derivation. A typical element of $\mathbb{R} \langle
u \rangle$ is a rational function of $u$, $\dot{u}$, \dots,
$u^{(\nu)}$, \dots, with real coefficients. 
\end{remark}
It is known (\cite{kolchin}) that the sum and the product of two elements which
are differentially algebraic over $\mathbb{R} \langle u \rangle$ is
again differentially algebraic over $\mathbb{R} \langle u \rangle$.
It is obvious moreover that any constant element, which satisfies
$\dot{y} = 0$, belongs to $\mathfrak{A}$.

Take two distinct points $(\tau,u), (\tau^\prime, u^\prime) \in
\mathcal{I} \times \mathcal{C}$. If $\tau \neq \tau^\prime$, then $y
= t$, which satisfies $\dot{y} = 1$, separates the two points. If
$\tau = \tau^\prime$, then assume that $u \neq u$ on the interval
$[0, \tau]$. It follows from Lerch's theorem (\cite{lerch}) (see,
also, \cite{miku}) that there exists a non-negative integer $\nu$
such that
$$
\int_{0}^{t} \sigma^\nu
u(\sigma) d \sigma \neq \int_{0}^{t} \sigma^\nu u^\prime(\sigma) d
\sigma
$$
The classic Cauchy formula demonstrates the existence of a
non-negative integer $\nu$ such that $y$, which satisfies $y^{(\nu)}
= u$, separates $(\tau,u)$, $(\tau, u^\prime)$.

This proof, which mimics to some extent
\cite{udine,smf} (see, also, \cite{suss}), shows that
$\mathfrak{A}$ is dense in $\mathfrak{S}$.

\section{Justification of the ultra-local model}
Assume that our SISO system is ``well'' approximated by a system
described by Equation \eqref{eq}. Let $\nu$ be a non-negative
integer such that
$$
\frac{\partial E}{\partial y^{(\nu)}} \not\equiv 0
$$
The implicit function theorem yields then locally
$$
y^{(\nu)} = \mathcal{E}(y, \dot{y}, \dots, y^{(\nu - 1)}, y^{(\nu +
1)}, \dots, y^{(a)}, u, \dot{u}, \dots, u^{(b)})
$$
It may be rewritten as Equation \eqref{ultralocal}.

\end{document}